\documentclass[11pt,a4paper]{article}
\pdfoutput=1
\usepackage[T1]{fontenc}
\usepackage[utf8]{inputenc}
\usepackage{lmodern}
\usepackage[english]{babel}
\usepackage[margin=2.4cm]{geometry}
\usepackage{microtype}
\usepackage{amsmath,amssymb,amsthm}
\usepackage{booktabs,array,tabularx}
\usepackage{enumitem}
\usepackage[authoryear,round]{natbib}
\usepackage{tikz}
\usetikzlibrary{calc,decorations.pathreplacing,arrows.meta,positioning}
\usepackage{xcolor}
\usepackage[hidelinks]{hyperref}
\usepackage{url}
\makeatletter\g@addto@macro\UrlBreaks{\do\-}\makeatother

\hypersetup{
  pdftitle={An event is a set, not a number: building the algebra of events before the probability measure in a school progression},
  pdfauthor={Vincent Tariel-Adam},
  pdfkeywords={probability, events, sample space, algebra of events, didactics of mathematics, curriculum design, modelling}
}

\newtheorem{principle}{Design principle}
\theoremstyle{definition}
\newtheorem{hypothesis}{Hypothesis}

\definecolor{evcol}{RGB}{110,150,200}   
\definecolor{numcol}{RGB}{225,150,90}   
\definecolor{greycol}{RGB}{150,150,150}
\tikzset{
  outcome/.style={circle,draw,inner sep=0pt,minimum size=4.2mm,font=\scriptsize},
  cell/.style={draw,minimum size=5.2mm,inner sep=0pt,font=\scriptsize},
}

\title{An event is a set, not a number:\\ building the algebra of events before the probability measure\\ in a school progression}
\author{Vincent Tariel-Adam\\[2pt]
  \small Independent researcher; CommeUnJeu (\url{https://commeunjeu.org}), France\\
  \small \texttt{vincent.tariel@alumni.polytechnique.org}}
\date{October 2026}

\begin{document}
\maketitle

\begin{abstract}
\noindent
Asked what an event is, pupils often answer with a fraction. This paper argues that school probability should separate two learnings and order them: the algebra of events (random experiment, outcomes, sample space, events as sets of outcomes, the operations not, or, and) and the probability measure, which attaches a number to each event. Before learning to compute the probability of an event, pupils should learn to build, identify and handle the event. The thesis rests on the structure of finite probability, in which the rules are facts about events carried over by an additive function, and on the history of the subject, in which errors such as d'Alembert's two tosses concerned the weighting of outcomes and were resolved by refining events into equally likely outcomes. Five design principles are stated: two objects with two names; the event before its probability; modal judgements (certain, possible, impossible) on sets before numbers; the algebra of events taught as such, with the rules derived from it; the sample space as a modelling choice, with frequencies last, as statistical estimates, and the separation repeated for conditional probability and random variables. They are confronted with the probability chapters, from age 9 to 18, of the French-curriculum courses published at \url{https://commeunjeu.org}, with research and other curricula. In primary school events precede probabilities without set notation, which comes at age 12; the lessons apply the principles and prove the rules from events in lower secondary school. But, in the version analysed, events account for 18 to 44 percent of the exercises, proofs are hidden in upper secondary school, and automatically corrected calculations record only the number. It proposes number-free tasks on events for ages 6 to 8 and changes to the International Baccalaureate curriculum. No empirical evaluation is reported; hypotheses are stated.
\end{abstract}


\section{Introduction}\label{sec:intro}

In the author's teaching experience, pupils of thirteen who have ``done probability'', asked what an event is, often answer with a number: ``an event is one half'', ``the event is three out of six''. The chapters analysed below list the same confusion among ``two classic mistakes''; how frequent it is remains an empirical question (Hypothesis~\ref{h:object}). A plausible reading is that such a pupil has met the event only inside a question---``what is the probability of rolling an even number?''---whose expected answer is a fraction, and that the two objects fused in that sentence have never been separated. One of them is a \emph{set}: rolling an even number is the set of outcomes $\{2,4,6\}$ of the experiment. The other is a \emph{number}: the probability $\frac{1}{2}$ that a model attaches to that set. The everyday way of speaking about chance tends to dissolve the object into its measure: ``the probability of an even number'' names the event only as the argument of a number. French didactic research met a similar difficulty with area, where a surface, its area as a magnitude and the number that measures it have to be told apart \citep{douady1989}; the event and its probability are the analogue of the surface and its measure.

This paper argues that school probability should separate two learnings and order them. The first is the \emph{algebra of events}: recognising a random experiment, listing its possible outcomes, forming the sample space, describing an event as a set of outcomes, and operating on events with the words \emph{not}, \emph{or}, \emph{and}---the complement, the union, the intersection---together with the two extreme events, the impossible and the certain. The second is the \emph{probability measure}: attaching to each event a number between $0$ and $1$ in a way that respects those operations. The thesis is that, before pupils learn to compute the probability of an event, they should learn to construct, identify and manipulate the event itself; and that this order is the condition of a genuine modelling approach to probability, rather than the learning of a ritual---``favourable over total''---that produces fractions about objects the pupil has never built.

The thesis is not new in mathematics: since \citet{kolmogorov1933}, events form an algebra of sets and probability is an additive function on it, and Kolmogorov illustrates the definition with two tosses of a coin, whose four outcomes $HH$, $HT$, $TH$, $TT$ make ``every event [\ldots] a set of elementary events'' \citep[\S2]{kolmogorov1950}. It is not new in research either. The review of \citet{bryant2012} calls working out the sample space ``the first and essential step in solving any probability problem'' and the lack of intervention studies on teaching it ``the most glaring gap'' in the field; an intervention study then showed that ten- and eleven-year-olds taught the sample space as an object of its own progressed more than comparison groups \citep{nunes2014}; and the framework of \citet{jones1997} treats the sample space and the probability of an event as distinct constructs. In France, the inter-IREM commission described an event as ``a particular set of possible outcomes'' \citep{henry2001experience}, and Girard noted that the programmes of the time had no sets at their disposal any more \citep{girard2001}. What this paper adds is an analysis of what the thesis requires of a progression over the whole of schooling, and a confrontation of a progression that claims to apply it with the chapters that implement it. The progression is that of the French-curriculum courses of the CommeUnJeu platform, from the fourth year of primary school to the final year of secondary school; the paper also examines the early years, which the French programmes leave without probability, and the International Baccalaureate (IB) curriculum, in which the separation is less developed.

\paragraph{Five design principles.} The progression rests on five principles, which the paper states, justifies and confronts with the chapters:
\begin{enumerate}[label=(\roman*)]
  \item two objects: an event is a set of outcomes, a probability is a number attached to it, and the two are never conflated in words or in notation (Section~\ref{sec:p1});
  \item the event first: in every chapter, outcomes and events precede probabilities, and every calculation begins by constructing the event---in primary school without any set formalism, the event being named by a sentence and built by listing the outcomes that make it happen, and from \emph{cinqui\`eme} in the language of sets (Section~\ref{sec:p2});
  \item the modal before the numerical: ``impossible'', ``certain'' and ``possible'' are first judged as relations between an event and the set of possible outcomes---none, all, some---before any number is attached (Section~\ref{sec:p3});
  \item the algebra of events is taught as such, and the rules of probability are carried over from it: the complement rule from the partition into an event and its complement, the addition relation from the double counting of the intersection (Section~\ref{sec:p4});
  \item the universe is a model: the choice of the outcomes and the hypothesis of equally likely outcomes are modelling decisions; the frequency of an event in a series of trials is a statistic, which estimates the probability or tests the model, and it comes at the end of the chapter, as a related notion that links probability to statistics; and the separation of the two objects is repeated at higher levels, where conditional probability changes the measure and not the events, and a random variable is a function before it has a law (Section~\ref{sec:p5}).
\end{enumerate}

\paragraph{Genre and scope.} The paper is a design rationale and a curriculum analysis. Like a hypothetical learning trajectory \citep{simon1995} or a research-based learning trajectory \citep{clements2004}, the progression embodies learning goals, an order of activities and hypotheses about learning; the paper makes them explicit, relates them to the mathematics, the history and the research literature, and confronts them with the materials. Its method borrows the \emph{a priori analysis} of didactic engineering \citep{artigue1990,artigue1992} at the level of global choices---which object comes first, which representation, which order---and gives one task-level analysis (Section~\ref{sec:cycle2}). It contains no classroom experiment: statements about learners are taken from the cited literature or stated as hypotheses (Section~\ref{sec:hypotheses}). The order it defends is not the deductive order of an axiomatic theory: events are built from situations---a bag, a die, two tosses---and frequencies, which some programmes put first, are frequencies \emph{of events}, met, from \emph{sixi\`eme} on, at the end of the chapter as statistical estimates of the probability of an event already built. None of the principles is new on its own; the contribution lies in their articulation over the whole of schooling, in the explicit statement of what each choice requires, and in an audit of the materials that shows where they carry the design and where they do not.

\paragraph{Corpus and method.} The corpus is the probability chapters of the French-curriculum courses of the platform, one per class from CM1 to the \emph{seconde} and several in the \emph{premi\`ere} and the \emph{terminale}, each with a lesson and an exercise set, in the version of 10~October~2026; the reservoir chapters of the earlier years from which they and the IB chapters are assembled; and the IB chapters of the same platform, from PYP~0 to the Diploma Programme. Unless a passage is explicitly identified as describing a later revision, descriptions, counts and quotations refer to this version; Sections~\ref{sec:audit} and~\ref{sec:ib} report the revisions made on 16~October~2026 in response to the analysis. The curriculum chapters are assembled from a reservoir of typed blocks (definition, proposition, proof, method, example, exercise, answer), each with a stable identifier; a curriculum chapter is a list of links to reservoir blocks, a link following a chain until a block with a body is reached, and a flag can hide the proof attached to a proposition. Statements about the corpus are of three kinds: \emph{presence} (a block with a given content is displayed in a class), \emph{order} (one block precedes another), and \emph{counts}. For the counts of Section~\ref{sec:audit}, each exercise is coded by the competence under which its \emph{root} block---the end of its chain of links---is filed in the reservoir, so that the same block receives the same code in every class. The categories are: \emph{events} (recognise an experiment, list outcomes, build a sample space, identify, form or combine events); \emph{modelling} (recognise equally likely outcomes, decide whether a game is fair by comparing winning outcomes, independence of repetitions); \emph{likelihood words} (place an event or a probability on a scale of words); \emph{numerical} (compute a probability, compare, convert, estimate by frequency); and \emph{problems}. The coding sheet, which lists every exercise with its root block, competence and category, is available with the source documents. Written exercises, corrected against a rubric, are also examined for whether the rubric assesses the construction of the event. Counts describe what the chapters ask, not the time pupils spend, and a numerical exercise may involve the construction of an event; they are made on the English half of the bilingual sources. The French programmes in force are the 2025 programmes of cycles~2 and~3 \citep{men2024cycle2,men2025cycle3} and the 2026 programme of cycle~4 \citep{men2026cycle4}; the IB documents are the Primary Years Programme scope and sequence, the Middle Years Programme guide and the Diploma Programme guides \citep{ibo2009pyp,ibo2020myp,ibo2019aa}.

\paragraph{Class names.} French classes are named by their French name followed, at first use, by the year of schooling counted from CP: CP (Grade~1, age 6), CE1 and CE2 (Grades~2--3), CM1 and CM2 (Grades~4--5), \emph{sixi\`eme} to \emph{troisi\`eme} (Grades~6--9, written 6e to 3e), \emph{seconde}, \emph{premi\`ere} and \emph{terminale} (Grades~10--12). The kindergarten year before CP is called GS (\emph{grande section}, age~5). IB years are named PYP~0--5, MYP~1--5 and DP~1--2.

Section~\ref{sec:maths} sets out the mathematical structure on which the thesis rests and Section~\ref{sec:history} the historical record; Section~\ref{sec:background} reviews the research on learners; Sections~\ref{sec:p1}--\ref{sec:p5} take the principles in turn; Section~\ref{sec:progression} describes the progression class by class; Section~\ref{sec:audit} confronts it with its exercises; Sections~\ref{sec:cycle2} and~\ref{sec:ib} examine the early years and the IB; Section~\ref{sec:misconceptions} lists the misconceptions that the design addresses, and Section~\ref{sec:hypotheses} states hypotheses for an empirical evaluation.

\section{The mathematical structure: two worlds and one function}\label{sec:maths}

\paragraph{The finite model.} School probability is almost entirely finite probability. A random experiment is modelled by a finite set $\Omega$, the \emph{sample space} or universe, whose elements are the \emph{outcomes}; an \emph{event} is a subset of $\Omega$, and the events form the set $\mathcal{P}(\Omega)$ of all subsets. A \emph{probability} on $\Omega$ is a function $P:\mathcal{P}(\Omega)\to[0,1]$ such that $P(\Omega)=1$ and $P(A\cup B)=P(A)+P(B)$ whenever $A\cap B=\emptyset$. In the finite case, such a function is determined by the probabilities $p_\omega=P(\{\omega\})$ of the singletons, non-negative and of sum $1$, through $P(A)=\sum_{\omega\in A}p_\omega$; when all the $p_\omega$ are equal, $P(A)=\frac{|A|}{|\Omega|}$, Laplace's ratio of favourable to possible cases. Nothing in the definition prevents an outcome from having probability $0$, in which case a non-empty event can have probability $0$ and an event other than $\Omega$ probability $1$; in the models that school probability chooses for an experiment every outcome has a positive probability, and only then do ``probability $0$'' and ``impossible'' coincide. The distinction between the impossible event $\emptyset$ and an event of probability $0$ is itself a distinction between the two worlds. These definitions are those of the \emph{seconde} and \emph{premi\`ere} chapters analysed below (Section~\ref{sec:progression}) and of the DP syllabus \citep{ibo2019aa}; they are the finite form of the axioms of \citet{kolmogorov1933}.

\begin{figure}[ht]
\centering
\begin{tikzpicture}[scale=0.92]
  \draw[evcol!80!black,thick,rounded corners=6pt] (-5.6,-2.0) rectangle (-0.4,2.1);
  \node[evcol!70!black,font=\small\bfseries] at (-3,1.75) {events: subsets of $\Omega$};
  \draw[evcol!80!black,thick] (-3,-0.15) ellipse (2.15 and 1.45);
  \node[evcol!80!black] at (-4.75,1.0) {$\Omega$};
  \fill[evcol!25] (-3.55,-0.15) ellipse (1.0 and 0.75);
  \draw[evcol!80!black] (-3.55,-0.15) ellipse (1.0 and 0.75);
  \foreach \x/\y/\l in {-3.9/0.15/2,-3.3/0.3/4,-3.6/-0.55/6}{\node[outcome,fill=white] at (\x,\y) {\l};}
  \foreach \x/\y/\l in {-1.95/0.55/1,-1.75/-0.3/3,-2.35/-1.05/5}{\node[outcome,fill=white] at (\x,\y) {\l};}
  \node[evcol!80!black,font=\small] at (-3.55,0.85) {$A$};
  \draw[numcol!80!black,thick,rounded corners=6pt] (1.4,-2.0) rectangle (5.8,2.1);
  \node[numcol!70!black,font=\small\bfseries] at (3.6,1.75) {numbers in $[0,1]$};
  \draw[-{Stealth[length=2mm]}] (1.9,-0.15) -- (5.4,-0.15);
  \foreach \x/\l in {2.1/0,5.0/1}{\draw (\x,-0.3) -- (\x,0.0) node[above] {\small $\l$};}
  \fill[numcol!90!black] (3.55,-0.15) circle (2.4pt);
  \node[numcol!90!black,below=3pt] at (3.55,-0.15) {\small $P(A)=\tfrac12$};
  \draw[-{Stealth[length=2.5mm]},very thick,numcol!80!black] (-2.45,-0.15) .. controls (0,0.9) and (1.5,0.6) .. node[above,pos=0.55] {$P$} (3.45,0.0);
\end{tikzpicture}
\caption{Two worlds and one function. For one roll of a die, the event ``an even number'' is the set $A=\{2,4,6\}$ of three outcomes of $\Omega=\{1,\dots,6\}$; the probability is the number $P(A)=\frac{1}{2}$ that the model attaches to it.}
\label{fig:worlds}
\end{figure}
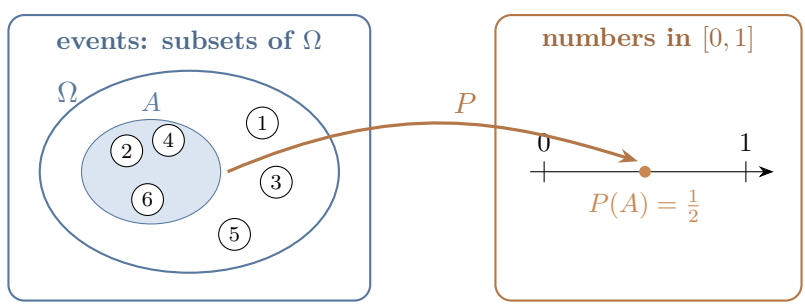

\paragraph{Two worlds.} The definitions describe two worlds joined by one function (Figure~\ref{fig:worlds}). In the world of events live $\Omega$ and its subsets, with their own operations---complement $\overline{A}$, union $A\cup B$, intersection $A\cap B$---their own extreme elements, $\emptyset$ (the impossible event) and $\Omega$ (the certain event), and their own facts, none of which mentions chance: $A$ and $\overline{A}$ partition $\Omega$; $A\cap B\subseteq A\subseteq A\cup B$; the events $A\cap B$, $A\cap\overline{B}$, $\overline{A}\cap B$ partition $A\cup B$. Under the dictionary \emph{not}/\emph{or}/\emph{and} $\leftrightarrow$ complement/union/intersection, this world is the Boolean algebra of the propositions that can be stated about the outcome of the experiment, which is how \citet{boole1854} first treated it. In the world of numbers live the real numbers of $[0,1]$. The probability $P$ belongs to neither world: it is the bridge that carries each event to a number.

\paragraph{The rules are facts about events, carried over.} Every rule of elementary probability is a fact of the world of events transported by $P$. Because $A$ and $\overline{A}$ are disjoint and cover $\Omega$, additivity gives $P(A)+P(\overline{A})=P(\Omega)=1$, the complement rule. Because $A\cup B$ is the disjoint union of $A$ and $B\cap\overline{A}$, and $B$ the disjoint union of $B\cap\overline{A}$ and $A\cap B$, additivity gives $P(A\cup B)+P(A\cap B)=P(A)+P(B)$, the addition relation that the French programme states in that symmetric form in \emph{troisi\`eme} \citep{men2026cycle4}; in the world of events it says that the overlap was counted twice. Because $A\cap B\subseteq A$, additivity gives $P(A\cap B)\le P(A)$: the conjunction of two events is never more probable than either of them, which is the rule that intuitive judgements violate in the conjunction fallacy \citep{tversky1983}. Because the events $E_1,\dots,E_n$ of a partition are disjoint and cover $\Omega$, $\sum_i P(E_i)=1$, and for any event $B$ the events $B\cap E_i$ are disjoint and cover $B$ (some may be empty), whence the law of total probability. In each case, the event-side fact and the axioms of $P$---additivity, together with non-negativity and $P(\Omega)=1$---give the rule; the design hypothesis is that a pupil who has the event-side fact can reconstruct the rule rather than memorise it.

\paragraph{The separation repeats.} Two later notions are the same separation one level up. Conditional probability given an event $A$ with $P(A)>0$ can be seen in two equivalent ways. In the first, $A$ becomes a new universe and the probabilities of its subsets are rescaled so that they add up to $1$; this is the ``new population'' of the two-way tables. In the second, $P_A(B)=P(A\cap B)/P(A)$ is a \emph{new probability} on the \emph{same} events of $\Omega$, which gives probability $0$ to the event $\overline{A}$, non-empty as soon as $A\neq\Omega$: knowing that $A$ has occurred changes the measure, not the world of events. The second view keeps $A$ and $B$ in one algebra, which is what Bayes' formula and the comparison of $P_A(B)$ with $P_B(A)$ require; it also shows that, once a measure is fixed, ``probability $0$'' and ``impossible'' can differ, which is one more reason to keep the two worlds apart. A random variable is a \emph{function} $X:\Omega\to\mathbb{R}$; the notation $(X=k)$ names an event, the set of outcomes at which $X$ takes the value $k$; and the law of $X$ is the old bridge applied to these events. Confusions between $P(A\cap B)$ and $P_A(B)$ are well documented \citep{batanero2016}; the design hypothesis is that keeping the two worlds apart makes each object of these formulas easier to place.

\paragraph{A model, not a fact.} The choice of $\Omega$ and of the $p_\omega$ is a modelling decision. The equality of the $p_\omega$, for a die or a coin, is a hypothesis justified by a symmetry and confronted with experience; for a drawing pin or a cone it is false, and the probabilities are estimated by frequencies. The French programme of cycle~4 says it in so many words: equiprobability ``is a hypothesis which cannot be proved, but which can be judged relevant by an argument of symmetry'', and ``in every case one distinguishes the real situation from the probability model used to describe it'' \citep{men2026cycle4}. French didactic research has made this distinction between the situation and its model the organising idea of the teaching of probability at secondary level \citep{henry2001}. The thesis of this paper is its consequence at the level of objects: the model has two parts, a universe with its events and a measure, and the first part must exist before the second can be chosen.

\paragraph{An analogy for teachers.} A shape is not its area. Shapes can be cut, glued and complemented; areas are numbers that respect those operations, the area of a disjoint union being the sum of the areas. Pupils do confuse a surface with its area and its measure, and teaching area means separating them \citep{douady1989}; the confusion of ``rolling an even number'' with $\frac{1}{2}$ is of the same kind. The analogy is the measure-theoretic one, domesticated: probability is the geometry of chance, events are the shapes and $P$ is the area. It also explains why Venn diagrams, in which an event is a region, are the natural representation of the algebra of events; a probability can be read on them as a share of area only when the regions are drawn with areas proportional to probabilities.

\section{The historical record: weighting outcomes, refining events}\label{sec:history}

The history of probability is often told as the history of a calculus. Read from the point of view of this paper, it also records a long difficulty that lies upstream of the calculus: several of the famous errors and disputes of the subject turn on which outcomes are taken to be equally likely, and their resolution consisted in refining coarse cases into finer ones---an operation on events---before any arithmetic. Not every classical difficulty is of this kind: the St Petersburg paradox concerns expectation, the disputes over inverse probability concern prior probabilities, and puzzles such as Monty Hall concern the protocol of conditioning. The cases below are those in which the world of events is at stake.

\paragraph{Three dice.} Galileo's short note on dice starts from an observation that he attributes to players: with three dice, a total of $10$ comes up more often than a total of $9$, although each total can be written in six ways as a sum of three numbers from $1$ to $6$ \citep{galilei1898,david1962}. Whether players could detect the difference empirically has been doubted; Galileo's explanation does not depend on it. He counts not the six partitions but the ordered triples: $27$ of the $216$ possible throws give $10$ and only $25$ give $9$ \citep[art.~8]{todhunter1865}. The six ``ways'' are events of unequal size in the model whose $216$ ordered throws are equally likely. Leibniz made the same kind of error with two dice, judging $11$ and $12$ equally easy because each ``can be made in only one way'' \citep[art.~77]{todhunter1865}; the same confusion between an unordered description and a set of ordered outcomes---``a five and a six'' against ``two sixes''---is still made by most pupils and by prospective teachers \citep{fischbein1997}.

\paragraph{The problem of points.} In 1654 Pascal and Fermat exchanged letters on the division of the stakes of an interrupted game \citep{fermat1894,devlin2008}. Fermat's method was to extend the game, in thought, to the number of rounds after which it would certainly be decided, to list all the sequences of results of those rounds, and to count those favourable to each player. Roberval objected to the fiction of rounds that would never be played, and Pascal at first thought that the method failed for three players. Fermat answered that ``everything that happens after one of the players has won is of no further use'', and that ``this fiction of extending the game to a certain number of rounds serves only to facilitate the rule and, in my opinion, to make all the chances equal'' \citep[letter of 25 September 1654]{fermat1894}. The fictitious extension is a choice of sample space made so that its outcomes are equally likely; each real end of the game is an event, a set of these outcomes. The arithmetic was never in dispute.

\begin{figure}[ht]
\centering
\begin{tikzpicture}
  \begin{scope}
    \node[font=\small] at (1.05,1.55) {second toss};
    \node[font=\small,rotate=90] at (-1.05,0.0) {first toss};
    \node[font=\small] at (0.45,1.05) {H}; \node[font=\small] at (1.65,1.05) {T};
    \node[font=\small] at (-0.55,0.45) {H}; \node[font=\small] at (-0.55,-0.75) {T};
    \foreach \x/\y/\l/\f in {0/0/HH/1,1.2/0/HT/1,0/-1.2/TH/1,1.2/-1.2/TT/0}{
      \ifnum\f=1 \fill[evcol!30] (\x,\y) rectangle ++(1.2,0.9); \fi
      \draw (\x,\y) rectangle ++(1.2,0.9); \node[font=\small] at (\x+0.6,\y+0.45) {\l};}
    \node[font=\small,align=center] at (1.2,-2.0) {four equally likely outcomes:\\ ``at least one head'' $=\{HH,HT,TH\}$, $P=\tfrac34$};
  \end{scope}
  \begin{scope}[xshift=6.3cm]
    \foreach \y/\l/\f in {0.5/{H (game over)}/1,-0.35/{T then H}/1,-1.2/{T then T}/0}{
      \ifnum\f=1 \fill[numcol!30] (0,\y) rectangle ++(2.8,0.7); \fi
      \draw (0,\y) rectangle ++(2.8,0.7); \node[font=\small] at (1.4,\y+0.35) {\l};}
    \node[font=\small,align=center] at (1.4,-2.0) {d'Alembert's three cases, counted\\ as equally likely: $P=\tfrac23$};
  \end{scope}
\end{tikzpicture}
\caption{Two tosses of a coin. Left: the sample space of four equally likely ordered pairs; the event ``at least one head'' is a set of three of them. Right: d'Alembert's three cases, a coarser sample space whose outcomes are not equally likely (the first gathers $HH$ and $HT$), counted as if they were.}
\label{fig:dalembert}
\end{figure}
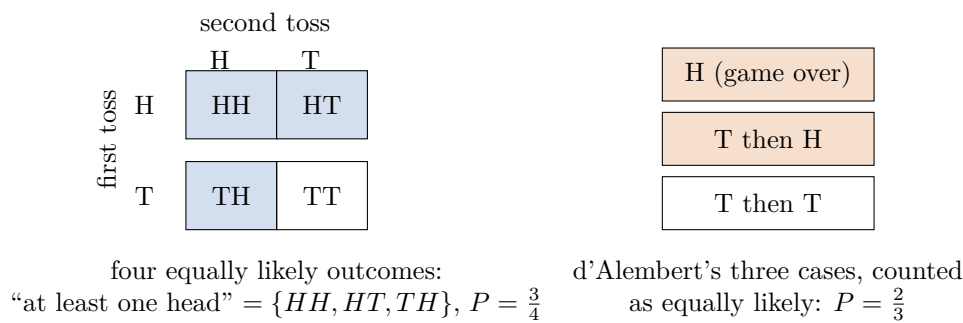

\paragraph{Two tosses of a coin.} In the article ``Croix ou pile'' of the \emph{Encyclop\'edie} (1754), d'Alembert asked for the probability of obtaining heads at least once in two tosses \citep{dalembert1754}. He argued that there were three cases---heads at the first toss, which ends the game; tails then heads; tails twice---of which two are favourable, and concluded $\frac{2}{3}$ instead of $\frac{3}{4}$ (Figure~\ref{fig:dalembert}). His three cases form a legitimate sample space; what is wrong is to treat them as equally likely \citep{henry2001croix}. Laplace's answer, sixty years later, is an operation on events: the first case ``is a simple event which corresponds to the two compound events'' heads--heads and heads--tails, and the extended game ``serves only to reduce the various cases to equally possible cases'' \citep[pp.~7--8]{laplace1814}. Weighted correctly, $\frac12$, $\frac14$, $\frac14$, the three cases give $\frac34$; the way to find these weights is to decompose the first case into two of the four equally likely ordered pairs.

\paragraph{Equally possible cases.} Laplace's classical definition, the ratio of the number of favourable cases to the number of all possible cases, ``equally possible'' \citep{laplace1814}, puts the difficulty into a single word: the definition presupposes that the cases have been chosen so that they are equally likely, which is a property of the measure, used to define the measure. The circularity was felt throughout the nineteenth century, and it is the same one that a pupil meets when told to count ``favourable over total'' without being told which outcomes to count.

\paragraph{Bertrand's chords.} \citet{bertrand1889} asked for the probability that a chord drawn ``at random'' in a circle is longer than the side of the inscribed equilateral triangle, and gave three reasoned answers, $\frac{1}{3}$, $\frac{1}{2}$ and $\frac{1}{4}$, according to how the chord is chosen; choosing ``at random among an infinity of possible cases'', he wrote, ``is not a sufficient indication''. The three answers correspond to three probability measures on the same set of chords, obtained from three ways of parametrising them. The paradox shows that ``at random'' does not define a probability until the outcomes and the measure on them are specified.

\paragraph{The axiomatisation.} When \citet{kolmogorov1933} gave probability its axioms, the first object he introduced was not a number but a set $E$ of elementary events and a field $\mathfrak{F}$ of its subsets, the random events; the probability came second, as a non-negative additive function on $\mathfrak{F}$ with value $1$ on $E$. His first illustration is d'Alembert's experiment, two tosses of a coin, and his third section opens on ``random events, as sets'', with a dictionary between the language of sets and the language of events \citep[\S\S2--3]{kolmogorov1950}---the dictionary that the French \emph{seconde} chapter reproduces (Section~\ref{sec:p1}). The construction answered a question that Hilbert had included among the problems of 1900, the axiomatisation of the probability calculus, and it gathered work of the preceding decades on measure and integration \citep{shafer2006}. Its effect on the old disputes was to make them statable: a disagreement about a probability becomes either a disagreement about the model---which set of outcomes, which measure---which the formalism makes explicit and open to argument and to experiment, or an error within a model, which it settles.

Read in this light, the history supports the thesis in a precise form. A correct sample space does not by itself prevent wrong weights---d'Alembert's three cases are a correct sample space---but the remedy that Galileo, Fermat and Laplace applied was in each case a refinement of coarse cases into finer, equally likely ones, which presupposes that one can decompose an event into outcomes and recombine them. The history is used in the progression as a source of problems in which this construction is the whole difficulty: Galileo's three dice, d'Alembert's two tosses, the problem of points and the Chevalier de M\'er\'e's reasoning appear in the \emph{troisi\`eme} chapter under the heading \emph{Solving Historical Problems of Chance}, and the French programme of \emph{troisi\`eme} itself proposes ``the problem of d'Alembert's error'' as a historical extension \citep{men2026cycle4}.

\section{What research says about learners}\label{sec:background}

\paragraph{Development and instruction.} The first systematic study of children's ideas of chance, by \citet{piaget1951}, tied the notion of chance to the development of operational thought, and in particular to the combinatorial operations on which the systematic listing of possible cases depends. \citet{fischbein1975} qualified the conclusion in a way that matters for curriculum design: children have \emph{primary intuitions} of chance well before the formal stage, and instruction can build \emph{secondary intuitions}; his book gives combinatorial analysis a chapter of its own. Fischbein and Schnarch later showed that, without instruction, some misconceptions grow \emph{stronger} with age and one is stable: asked whether ``a five and a six'' or ``two sixes'' is more likely with two dice, $70$ to $78\%$ of pupils from the fifth to the eleventh grade, and of prospective mathematics teachers, answered that they are equally likely \citep{fischbein1997}---Leibniz's error, an unordered description taken for a single outcome. Green's survey of nearly three thousand English pupils aged 11 to 16 gave the first large-scale picture of pupils' probability concepts at secondary level \citep{green1983}.

\paragraph{The sample space as an object of teaching.} \citet{bryant2012} identify four cognitive demands of probability for children---understanding randomness, working out the sample space, comparing and quantifying probabilities, understanding correlations---and call the second ``the first and essential step in solving any probability problem'', while noting that it has been ``relatively neglected''. Their review also discusses \emph{aggregation}, the grouping of equally likely outcomes into events of unequal probability, such as the eleven totals of two dice made of thirty-six ordered pairs; in its vocabulary, the elements of the sample space are themselves called ``events'', which shows how easily the two levels merge even in research writing. \citet{nunes2014} proposed classification, logical multiplication and ratio as the conceptual bases of the sample space and showed, in an intervention with children aged ten and eleven, that those taught the sample space as an object of its own progressed more than a problem-solving group and an unseen comparison group, two months later included. The framework of \citet{jones1997} describes four constructs---sample space, probability of an event, probability comparisons, conditional probability---each with four levels of thinking, from subjective to numerical. Studies of the total of two dice and two spinners \citep{pratt2000} and of the ``combinations tower'' \citep{abrahamson2009} show learners reasoning about aggregate events while the designer reasons about the outcomes that compose them.

\paragraph{Young children.} Studies of children aged four to eight bear directly on Principle~\ref{pr:modal}. \citet{way2003} found that young children knew that a colour must be present in a bag for drawing it to be possible, but not that ``for an outcome to be impossible there must be none of that colour present at all'', and that, at the earliest stage, no clear connection was made between the structure of the sample space and the likelihood of particular events. \citet{nikiforidou2019} found that children of four to six judged likelihood better with tangible materials and ``when the sample space is simpler'', and a teaching experiment with five-year-olds identified the interpretation of ``different'' outcomes of a two-stage experiment as one of three critical steps \citep{kafoussi2004}.

\paragraph{Misconceptions and the world of events.} Several of the best-documented misconceptions can be read, at least in part, as difficulties with events or with the weighting of outcomes rather than with numbers. In the \emph{outcome approach} described by \citet{konold1989}, a question about a probability is read as a request to predict the outcome of the next trial. In the \emph{equiprobability bias} studied by \citet{lecoutre1992}, results are judged equally likely ``because it is a matter of chance''; Lecoutre could induce a better model by masking the chance aspect of the task, but transfer back to a random situation was weak. The \emph{representativeness} heuristic \citep{kahneman1972} judges a single outcome by its resemblance to the process instead of by the size of the event it belongs to. In the \emph{conjunction fallacy}, adults judge $A\cap B$ more probable than $A$; \citet{tversky1983} called the rule violated \emph{extensional}, because it holds ``because the extension [\ldots] of the conjunction is included in the extension of its constituents''. \citet{chernoff2009,chernoff2011} showed that judgements about sequences of coin tosses which look wrong can be correct statements about a different partition of the sample space---``mixed'' against ``not mixed'', for instance---so that the error lies in the alignment of the event described with the event intended. Conversely, presenting statistical information as natural frequencies---counts of individuals, that is, sizes of sets---improves Bayesian reasoning without instruction \citep{gigerenzer1995}. Several of the remedies proposed in these studies consist in making the events, or the sample space, visible.

\paragraph{Modelling, and an objection from within the field.} French didactic research has insisted on the distinction between a random situation and the probability model chosen to describe it \citep{henry2001,coutinho2001}. Henry dates a turn in the French programmes to 1991, when an a priori universe with Laplace's equiprobability gave way to an entry through observed frequencies \citep{henry2001experience}; Girard noted the cost of that turn, the loss of the sets on which events are built \citep{girard2001}. \citet{batanero2005nature} argue that the classical, frequentist and subjective meanings of probability are complementary and should all inform the curriculum. In the theory of didactical situations, the experiment of Brousseau and colleagues had pupils infer the hidden content of an opaque bottle from repeated observations \citep{brousseau2001}: the object of the model was the set of possible results before it was a frequency. \citet{steinbring1991} gives the strongest objection to an order of objects: probability is a theoretical concept that lives in the relation between a model and experiments, and ``there cannot be a logical and deductive course for introducing the basic concepts and then constructing the theory upon them''. The design accepts the premise and rejects the inference. Events are not introduced by definition before experience; they are built in situations---a bag, a die, two tosses---and the frequencies that the pupils compute, at the end of the chapter, are frequencies of events that they have named. What the design orders is not theory before experience but the object before its measure, within each situation. In the same spirit, \citet{freudenthal1973} warned against the ``anti-didactical inversion'' that teaches the finished structure of a theory instead of the activity that produces it; the tasks of Section~\ref{sec:cycle2} are designed as such an activity.

\paragraph{The objection of the ``modern mathematics''.} A proposal to teach events as sets from the primary years meets an objection that has a history. In France, the Lichnerowicz commission, created in 1966, prepared ``modern'' programmes built on sets and structures; the realities of mass secondary schooling produced an open crisis as early as 1971, and the reform was abandoned within a decade \citep{denfert2011,gispert2014}. In the United States, the ``new math'' met a similar fate \citep{kline1973}, and Freudenthal criticised the ``false set theory'' of school textbooks \citep{freudenthal1973sets}. The objection deserves a precise answer, because the hoops and diagrams of Section~\ref{sec:cycle2} resemble the activities of that period. Four features of the reform are generally held responsible for its failure: sets were taught as the foundation of all mathematics, before and independently of any problem; vocabulary and symbols were taught and assessed for their own sake; the activities were disconnected from situations in which they made a difference; and teachers were asked to teach a mathematics they had not been trained in. The proposal avoids each. The sets are always sets of possible results of one concrete experiment, and they answer a question about chance (will it happen? which is more likely?); a draw can refute a judgement, and inspecting the whole bag settles it, so that the situation, not the teacher, gives feedback---the \emph{milieu} of the theory of didactical situations; no symbol is used in the early years, and set notation comes only where the chapters introduce it, in \emph{cinqui\`eme}, one year before the French programme \citep{men2026cycle4}; and what is assessed is the judgement about chance, not the set vocabulary. The French programme of \emph{sixi\`eme} states that pupils are not expected to use the specific vocabulary---experiment, outcome, universe, event---on their own \citep{men2025cycle3}; the corpus uses this vocabulary from CM1, which is a deliberate departure that the hypotheses of Section~\ref{sec:hypotheses} put to the test.

\section{Principle 1: two objects, two names}\label{sec:p1}

\begin{principle}[An event is a set, a probability is a number]\label{pr:two}
The event and its probability are two objects of different natures---a set of outcomes and a number---and they are never conflated, in words or in notation.
\end{principle}

The principle is stated in the chapters as a rule of speech and writing. The first definition of an event, from CM1 on, is set-theoretic: ``An \emph{event} is a set of outcomes from the set of all possible outcomes.'' From \emph{cinqui\`eme} it becomes ``a \emph{subset} of the sample space'', denoted by a capital letter. The \emph{premi\`ere} chapter states the principle itself, in a text block titled \emph{An Event Is a Set and a Probability Is a Number}:
\begin{quote}\small
Rolling a fair die, the event ``the result is even'' is the set $A=\{2,4,6\}$; its probability is the number $P(A)=\frac{3}{6}=\frac{1}{2}$. Events combined with ``or'', ``and'', ``not'' measure nothing by themselves; a probability attaches a number to each event. The two objects must be kept apart, because two classic mistakes come from mixing them: writing $P(A)=\{2,4,6\}$ [\ldots]; writing $A=\frac{1}{2}$ or saying ``the event $A$ is one half'' [\ldots] The letter $P$ is a function: it takes an event $A\subseteq U$ and returns a number, its probability $P(A)$.
\end{quote}
The same care governs the probability of a single outcome. In \emph{quatri\`eme}, a text block explains that an outcome $x$ and its singleton event $\{x\}$ ``describe the same result: speaking of the chance of an outcome is therefore speaking of the chance of an event'', so that the probability of the outcome is $P(\{x\})$, the shorthand $P(x)$ being a permitted abuse. And in \emph{seconde}, a table sets the two vocabularies side by side, notation by notation: $\Omega$, ``universal set'', ``sample space''; $\emptyset$, ``empty set'', ``impossible event''; $E$, ``subset of $\Omega$'', ``event''; $E\cup F$, ``union'', ``$E$ or $F$''; $E\cap F=\emptyset$, ``disjoint'', ``mutually exclusive''.

Two features of the early classes deserve comment. In CM1 and CM2 the notation is $P(\text{event})$, with a phrase in place of the event---$P(\text{``choosing a red ball''})$. This is a deliberate choice: in primary school the event is a sentence together with the outcomes that make it true, without set notation, and the letter naming a set comes in \emph{cinqui\`eme}. And the line on which likelihood words are placed (\emph{impossible}, \emph{less likely}, \emph{even chance}, \emph{more likely}, \emph{certain}) is illustrated, in the shared CM1--\emph{cinqui\`eme} block, by two events of experiments (``rolling a die and getting a 3'', ``tossing a coin and getting heads'') and by three everyday happenings---``riding a dinosaur'', ``the sun will rise tomorrow''---that are not presented as events of an experiment whose outcomes are listed. A model exists for each of the latter (two outcomes, ``it happens'' and ``it does not''), but the chapter does not state it, so that the definition of an event given a few lines earlier is not applied to them. Section~\ref{sec:audit} returns to this point.

\section{Principle 2: the event first}\label{sec:p2}

\begin{principle}[The event before its probability]\label{pr:first}
In every chapter, the outcomes, the sample space and the events are introduced before any probability; and every calculation of a probability begins by constructing the event. In primary school the order holds without any set formalism: an event is named by a sentence and built by listing the outcomes that make it happen; the language and the notation of sets come in \emph{cinqui\`eme}.
\end{principle}

\paragraph{The order of the chapters.} The principle is visible in the section titles, which can be read as a statement of the thesis. CM1: \emph{Outcomes}, \emph{Events}, then \emph{Using Words to Describe Probability}, \emph{Using Numbers to Quantify Probability}, \emph{Calculating Probabilities}, \emph{Fair and Unfair Games}. CM2: \emph{Outcomes, Events and Probability}, \emph{Two-Step Experiments}, \emph{The Die Has No Memory}, \emph{Calculating and Comparing in Two Steps}. \emph{Sixi\`eme}: \emph{The Vocabulary of Probability} (outcomes, events, two-step experiments) before \emph{The Probability of an Event}. \emph{Cinqui\`eme} and \emph{quatri\`eme}: a first part titled \emph{The Algebra of Events}, then a part titled \emph{Probability}. \emph{Troisi\`eme} and \emph{seconde}: \emph{The Language of Events}, then \emph{Computing Probabilities} or \emph{The Law of Probability}. \emph{Premi\`ere}: \emph{Algebra of Events}, then \emph{Axioms and Rules of Probability}. In none of these chapters does a number attached to an event precede the definition of the event.

\paragraph{The order of a calculation.} The same order is prescribed for each calculation. The \emph{quatri\`eme} method \emph{Computing the Probability of an Event} has two steps:
\begin{quote}\small
1. \emph{Determine the event}: describe it, then list its outcomes---write it as a set. 2. \emph{Compute its probability}: if the outcomes are equally likely, $P(A)=\frac{\text{number of outcomes of } A}{\text{total number of outcomes}}$; otherwise, add the probabilities of the outcomes of the event; when the complementary event is simpler, $P(A)=1-P(\overline{A})$.
\end{quote}
The first step is the algebra of events, the second the measure; the three ways of computing differ only in the second step.

\paragraph{What the model answers do.} In the primary classes the worked answers apply the order without set notation, as the design intends. The CM1 lesson answers the example on even numbers by listing: ``The outcomes for `rolling an even number' are $2$, $4$ and $6$'', and its first calculation begins ``The event is `rolling an even number'. The outcomes that are even are $2$, $4$ and $6$. There are 3 favorable outcomes.'' A typical CM1 exercise answer reads: ``Each ball is equally likely to be drawn: count the balls, not the colours. The total number of possible outcomes [\ldots] is 5 [\ldots]. The number of outcomes for the event `choosing a red ball' is 2, as there are two red balls in the bag. Therefore, the probability [\ldots]''. The event is named and its outcomes are counted; the instruction ``count the balls, not the colours'' is a lesson about the sample space---the outcomes of the draw are the balls, not the colours---which is exactly the point at which the equiprobability bias enters. Set notation appears, as planned, from \emph{cinqui\`eme}; in the version analysed, the answer that the online calculation exercises recorded was, at every level, the number alone (Section~\ref{sec:audit}).

\section{Principle 3: the modal before the numerical}\label{sec:p3}

\begin{principle}[Certain, impossible, possible as relations between sets]\label{pr:modal}
The first judgements about chance are modal---certain, possible, impossible---and they are made by comparing the event with the set of possible outcomes: an event is certain when it contains all of them, impossible when it contains none, possible otherwise. Numbers come after.
\end{principle}

\begin{figure}[ht]
\centering
\begin{tikzpicture}
  \foreach \i/\cols/\word/\rel in {0/{red,red,red,red}/certain/{all are red: $E=\Omega$},1/{red,red,evcol,evcol}/maybe/{some are red: $\emptyset\neq E\neq\Omega$},2/{evcol,evcol,yellow!80!orange,evcol}/impossible/{none is red: $E=\emptyset$}}{
    \begin{scope}[xshift=\i*4.9cm]
      \draw[thick,rounded corners=8pt] (0,0) -- (0,1.5) .. controls (0.2,1.9) and (2.6,1.9) .. (2.8,1.5) -- (2.8,0) -- cycle;
      \foreach \c [count=\k] in \cols {\fill[\c!70] (0.1+0.55*\k,0.6) circle (0.2); \draw (0.1+0.55*\k,0.6) circle (0.2);}
      \node[font=\small\bfseries] at (1.4,-0.4) {\word};
      \node[font=\scriptsize,align=center] at (1.4,-0.9) {\rel};
    \end{scope}}
\end{tikzpicture}
\caption{The kindergarten method \emph{Picking without Looking}. The event is ``picking a red counter''; the outcomes are the counters. Certain, maybe and impossible are decided by comparing the red counters with all the counters---\emph{all}, \emph{some}, \emph{none}---before any number.}
\label{fig:bag}
\end{figure}
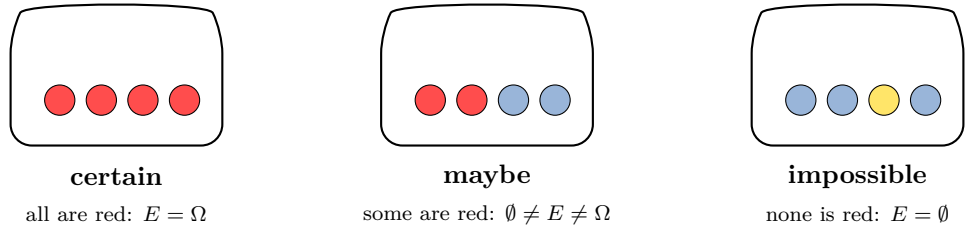

The clearest realisation of the principle is the earliest. The kindergarten chapter (GS, reservoir level used by the IB PYP~0 chapter) opens, like the later ones, with everyday happenings---the sun rising, rain tomorrow, a cat that flies---and then, in a section titled \emph{Picking without Looking}, gives a method, \emph{Looking into the Bag}, that is entirely set-theoretic without a single symbol:
\begin{quote}\small
We pick one counter from the bag \emph{without looking}. Before picking, we look at all the counters in the bag: if \emph{all} the counters are red, picking a red counter is \emph{certain}; if \emph{some} counters are red and some are not, picking a red counter: \emph{maybe}; if \emph{no} counter is red, picking a red counter is \emph{impossible}.
\end{quote}
The three quantifiers \emph{all}, \emph{some}, \emph{no} are the three possible relations between the event ``a red counter'' and the set of counters in the bag: equality with the whole, proper non-empty part, empty part (Figure~\ref{fig:bag}). The worked example adds the decisive remark: in a bag of two blue counters and one yellow, ``picking a red counter is impossible: there is no red counter in the bag''---an event of the experiment that contains no outcome. In \emph{seconde} the dictionary will say ``$\emptyset$, empty set, impossible event''; at age five the pupil decides impossibility by inspecting a set. The method answers the difficulty observed by \citet{way2003}: young children know that a colour must be present for drawing it to be possible, but not that it must be absent for drawing it to be impossible; the three quantifiers make the second condition as explicit as the first. In the version analysed, the IB PYP~1 and PYP~2 chapters reused the method, and a PYP~2 task asked the pupil to ``colour the marbles so that picking a red one is \emph{certain}''---the construction of a bag in which the event equals the sample space.

The later chapters keep the modal tier but change its referent. The probability line of CM1 to \emph{cinqui\`eme} orders likelihood words on a scale from impossible to certain and defines \emph{possible} as ``an event that may happen but is not sure to happen''; the definition of probability then places the number $0$ at ``impossible'' and $1$ at ``certain''. The identification of $0$ with ``impossible'' and $1$ with ``certain'' holds in the models that the chapters choose for an experiment, in which every outcome has a positive probability (Section~\ref{sec:maths}). The progression thus moves from a modal judgement made on a set (GS), to a scale of words, to a scale of numbers. What the line's illustrations lose, as noted above, is the set: ``riding a dinosaur'' is impossible, but no experiment is stated in which it is the empty event.

\section{Principle 4: the algebra of events, and the rules carried over}\label{sec:p4}

\begin{principle}[Operations on events first; rules of probability as their consequences]\label{pr:algebra}
The complement, the union, the intersection and the partition of the sample space are taught as operations on sets of outcomes, with the words \emph{not}, \emph{or}, \emph{and}; the rules of probability are presented as what additivity makes of facts about events.
\end{principle}

\paragraph{The operations.} The complementary event is defined in \emph{cinqui\`eme} as ``the set of all outcomes in the sample space that are \emph{not} in $E$'', with the notations $E'$ and $\overline{E}$. Union and intersection follow in \emph{quatri\`eme}, each with its word: $E\cup F$ is ``the event that occurs if $E$ occurs, or $F$ occurs, or both'', read ``$E$ or $F$'', the inclusive \emph{or} being stated; $E\cap F$ is read ``$E$ and $F$''. Mutually exclusive events are those whose intersection is empty, and a remark warns that ``mutually exclusive does not mean complementary'': with a die, $\{1,2\}$ and $\{5,6\}$ have no outcome in common but do not cover the sample space. These definitions are those of the French programme of \emph{quatri\`eme}, which asks pupils to ``use the vocabulary and the notations of sets to describe a random experiment [\ldots] and define the notion of event'' and to ``define complement, union, intersection, empty set (impossible event)'' \citep{men2026cycle4}. In \emph{premi\`ere} the partition of the sample space is introduced with a sentence that states the principle in its general form:
\begin{quote}\small
Cutting up $U$ is a question about \emph{sets of outcomes only}---it uses nothing about chance. Attaching a probability $P(E_i)$ to each piece is a \emph{separate} step---and \emph{however} we cut $U$ into pieces, the chances of the pieces always total $1$.
\end{quote}

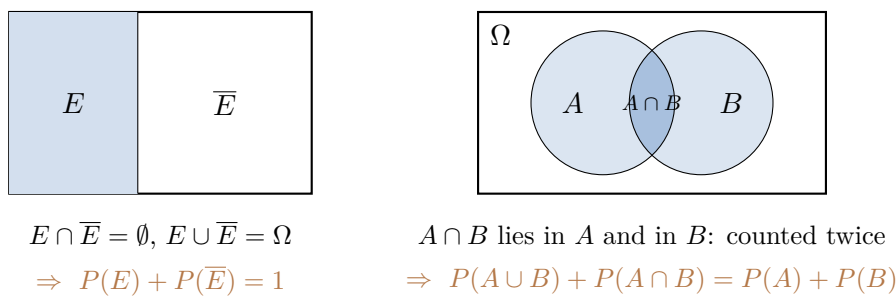
\begin{figure}[ht]
\centering
\begin{tikzpicture}
  \begin{scope}
    \draw[thick] (0,0) rectangle (4,2.4); \node at (0.3,2.1) {$\Omega$};
    \fill[evcol!30] (0,0) rectangle (1.7,2.4); \draw (1.7,0) -- (1.7,2.4);
    \node at (0.85,1.2) {$E$}; \node at (2.85,1.2) {$\overline{E}$};
    \node[font=\small,align=center] at (2,-0.55) {$E\cap\overline{E}=\emptyset$, $E\cup\overline{E}=\Omega$};
    \node[font=\small,align=center,numcol!80!black] at (2,-1.15) {$\Rightarrow\ P(E)+P(\overline{E})=1$};
  \end{scope}
  \begin{scope}[xshift=6.2cm]
    \draw[thick] (0,0) rectangle (4.6,2.4); \node at (0.3,2.1) {$\Omega$};
    \fill[evcol!25] (1.65,1.2) circle (0.95); \fill[evcol!25] (2.95,1.2) circle (0.95);
    \begin{scope} \clip (1.65,1.2) circle (0.95); \fill[evcol!60] (2.95,1.2) circle (0.95); \end{scope}
    \draw (1.65,1.2) circle (0.95); \draw (2.95,1.2) circle (0.95);
    \node at (1.25,1.2) {$A$}; \node at (3.35,1.2) {$B$}; \node[font=\scriptsize] at (2.3,1.2) {$A\cap B$};
    \node[font=\small,align=center] at (2.3,-0.55) {$A\cap B$ lies in $A$ and in $B$: counted twice};
    \node[font=\small,align=center,numcol!80!black] at (2.3,-1.15) {$\Rightarrow\ P(A\cup B)+P(A\cap B)=P(A)+P(B)$};
  \end{scope}
\end{tikzpicture}
\caption{Rules of probability as facts about events carried over by additivity. Left: an event and its complement partition the sample space. Right: in $A$ and $B$ taken separately, the outcomes of $A\cap B$ are counted twice.}
\label{fig:transport}
\end{figure}

\paragraph{The rules.} The intended relation between operations and rules is the one of Section~\ref{sec:maths}: each rule is an event-side fact transported by additivity (Figure~\ref{fig:transport}). The chapters of the \emph{coll\`ege} display it. In \emph{cinqui\`eme}, \emph{quatri\`eme} and \emph{troisi\`eme}, the complement rule is followed by its proof: ``every outcome of $U$ lies in \emph{exactly one} of $E$ and $E'$ [\ldots] The probability of an event is the sum of the probabilities of its outcomes. Adding $P(E)$ and $P(E')$ therefore counts the probability of every outcome exactly once''. In \emph{troisi\`eme}, the addition relation is followed by a Venn-diagram proof: ``if we simply add $P(E)$ and $P(F)$, we are double-counting the intersection region''. The chapters of the lyc\'ee hide these derivations: in \emph{seconde}, the complement rule and the addition relation are displayed without their proofs, and in \emph{premi\`ere} the complement rule and the sum over a partition are too, although the \emph{premi\`ere} chapter does display, with its proof, the monotonicity of probability, $A\subseteq B\Rightarrow P(A)\le P(B)$, which is the conjunction rule in set form. The \emph{troisi\`eme} and \emph{seconde} method for the addition relation also states the check that ``$P(A\cap B)$ can exceed neither $P(A)$ nor $P(B)$'', without the inclusion that justifies it. The rules are thus carried over from the events where they are first met, and stated as facts where they are used most.

\section{Principle 5: the universe is a model, and the separation repeats}\label{sec:p5}

\begin{principle}[The sample space as a modelling choice]\label{pr:model}
The outcomes are chosen: their choice and the hypothesis that they are equally likely are modelling decisions. The frequency of an event in a series of trials is a statistic: it estimates the probability of the event or tests the model, and it comes at the end of the chapter, once the event has been built and its probability defined in a model, known or not, as a related notion that links probability to statistics. The separation of events and measure is maintained at the higher levels, where conditional probability changes the measure on the same events and a random variable is a function before it has a law.
\end{principle}

\paragraph{Choosing the outcomes.} The progression meets the choice of outcomes three times before the lyc\'ee. In CM1, a remark compares two spinners: on the first, cut into six equal parts, two red parts give $P(\text{red})=\frac{2}{6}$; on the second, ``three colours, but $P(\text{red})$ is not $\frac13$: red covers half of the spinner'', and the probability is $\frac12$. The colours could serve as outcomes, but not as equally likely ones: the equally likely outcomes are the equal parts of the circle, and a colour names an event made of some of them---the same lesson as ``count the balls, not the colours''. Choosing outcomes that are equally likely, or weighting coarser outcomes correctly, is the modelling decision at stake. In CM2, the two-step experiment introduces outcomes that are ordered pairs, ``a result of the first step \emph{followed by} a result of the second'', listed by a table or a tree, and the method states the condition under which the table is complete: ``this works because the second step has the same outcomes whatever the first one gave''. In \emph{quatri\`eme}, the sample space of a two-step experiment is defined as ``the set of all these pairs'', with the qualification that if ``an earlier result rules some out, only the remaining pairs form the sample space''.

\paragraph{The hypothesis and the experiment.} Equiprobability is presented, from CM1 on, as a condition to check---``a fair die, a balanced spinner with equal sectors, marbles alike apart from their colour''---and its failure is met in \emph{cinqui\`eme} with a cone that lands on its base or its side, whose probability is estimated by frequency. In \emph{terminale}, a text block titled \emph{A Model, Not a Theorem} says of the independence of repeated trials that ``in a given situation it is a modelling choice, stated rather than deduced from the data''.

\paragraph{Frequencies at the end of the chapter.} From \emph{sixi\`eme} to \emph{seconde}, the last section or subsection of every chapter, before its exercises, deals with repeated trials: \emph{Repeating an Experiment} in \emph{sixi\`eme}, \emph{Frequency and Probability} in \emph{cinqui\`eme}, \emph{Frequencies and Probabilities: Fluctuation} in \emph{quatri\`eme}, \emph{From Chance to Experiment}---fluctuation, then simulation---in \emph{troisi\`eme}, and \emph{From the Model to the Sample}---frequencies, samples, sampling fluctuation, the law of large numbers---in \emph{seconde}; the CM1 and CM2 chapters have no frequencies. The place is a design decision, which follows from the thesis; its benefit over an earlier experimental entry is a hypothesis (Hypothesis~\ref{h:frequency}). The frequency of an event in $n$ trials is not a second definition of its probability but a statistic computed on a sample, and it serves two purposes, both of which presuppose the event and the model. It estimates a probability, in particular one that no symmetry gives: the cone lands on its base $15$ times in $50$, so that $P(\text{base down})\approx\frac{15}{50}$, and the lesson defines the ``experimental probability'' as ``an estimate found by repeating an experiment many times''. And it confronts a model with experience: in \emph{sixi\`eme}, a fair die gives no six in twelve rolls, and the method warns: ``A gap, on its own, does not prove that the die is loaded.'' The frequency is thus a related notion rather than a separate way into probability: it links probability to statistics, and the \emph{seconde} chapter makes the link explicit when it ends on the law of large numbers, where, for an individual drawn at random with equal chances, the proportion of a characteristic in the population equals the probability of drawing an individual who has it, ``so estimating one estimates the other''. The order of the chapter---the event, its probability in a model, then the frequency of the same event in a series of trials---carries the order of the thesis one step further. \citet{batanero2005nature} argue that the classical, frequentist and subjective meanings of probability are complementary; the progression keeps the frequentist meaning, as an estimate of a probability rather than as its definition. It also reverses the entry through observed frequencies of the 1991 French programme \citep{henry2001experience} without giving up what that entry sought, the confrontation of the model with experience: what the pupil confronts with the model is the frequency of an event already built. Simulation (\emph{troisi\`eme}), which repeats a model in place of the experiment, and the law of large numbers (\emph{seconde}, then \emph{terminale}), which justifies the estimate when the trials are independent and identical, belong to the same end of the chapter.

\begin{figure}[ht]
\centering
\begin{tikzpicture}
  \begin{scope}
    \draw[thick] (0,0) rectangle (4,2.4); \node at (0.3,2.1) {$\Omega$};
    \draw (1.6,1.2) ellipse (1.2 and 0.9); \draw (2.6,1.2) ellipse (1.1 and 0.8);
    \node at (0.9,1.2) {$A$}; \node at (3.3,1.2) {$B$};
    \node[font=\small,align=center] at (2,-0.5) {$P$ on all of $\Omega$};
  \end{scope}
  \draw[-{Stealth[length=2.5mm]},thick] (4.4,1.2) -- node[above,font=\small] {knowing $A$} (5.8,1.2);
  \begin{scope}[xshift=6.2cm]
    \fill[greycol!15] (0,0) rectangle (4,2.4); \draw[thick,greycol] (0,0) rectangle (4,2.4);
    \fill[white] (1.6,1.2) ellipse (1.2 and 0.9); \draw[thick] (1.6,1.2) ellipse (1.2 and 0.9);
    \begin{scope} \clip (1.6,1.2) ellipse (1.2 and 0.9); \fill[evcol!40] (2.6,1.2) ellipse (1.1 and 0.8); \end{scope}
    \draw[dashed] (2.6,1.2) ellipse (1.1 and 0.8);
    \node at (0.9,1.2) {$A$}; \node[font=\scriptsize] at (2.15,1.2) {$A\cap B$};
    \node[font=\small,align=center] at (2,-0.5) {$P_A(B)=P(A\cap B)/P(A)$:\\ same events, new measure};
  \end{scope}
\end{tikzpicture}
\caption{Conditional probability in the second view: the events are kept and the measure changes. The mass outside $A$ is set to zero and the rest is rescaled, so that $P_A(A)=1$ and the event $\overline{A}$, non-empty when $A\neq\Omega$, has probability $0$.}
\label{fig:conditional}
\end{figure}
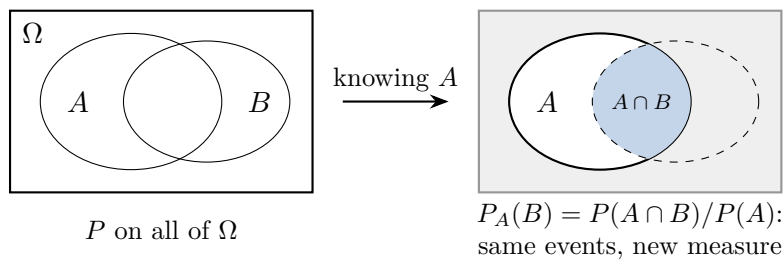

\paragraph{The separation one level up.} Conditional probability is introduced in \emph{seconde} through two-way tables, where ``the known event $B$ becomes the new population'', with the definition $P_E(F)=P(E\cap F)/P(E)$ and both of its notations, $P_E(F)$ and $P(F\mid E)$; the \emph{seconde} chapter already has weighted trees, the product rule on a path and the inversion of the conditioning, with screening tests, and \emph{premi\`ere} adds total probability and independence. The chapters take the first of the two views of Section~\ref{sec:maths}, the known event as a new universe, which is the natural one with counts and natural frequencies. The second view, a new probability on the same events, was not stated in the version analysed; it completes the principle where it is most needed, for Bayes' formula and for the distinction between $P_A(B)$ and $P(A\cap B)$ (Figure~\ref{fig:conditional}), and it has since been added to the \emph{premi\`ere} chapter (Section~\ref{sec:audit}). For random variables the principle is fully stated. The \emph{premi\`ere} chapter opens with a block titled \emph{From Events to Numbers}: ``An event is a set of outcomes and its probability a number. Many questions about a random experiment are about a number attached to each outcome rather than about an event''. The random variable is then defined as a function $X$ on the sample space; $(X=x)$ is defined as ``the set of outcomes where $X$ takes the value $x$''; and a text block explains that, once the variable is defined, ``one often sets the outcomes aside and keeps only what $X$ produces [\ldots] The outcomes come back whenever such a probability has to be computed.'' In \emph{terminale}, the sequence of $n$ trials has the Cartesian product of the sets of outcomes as its sample space, and the binomial coefficient counts ``the paths of the tree with $k$ successes'', that is, the outcomes of the event $(X=k)$.

\section{The progression class by class}\label{sec:progression}

Table~\ref{tab:progression} summarises the progression; the paragraphs below comment on its turning points.

\begin{table}[ht]
\centering\footnotesize
\begin{tabularx}{\textwidth}{>{\raggedright\arraybackslash}p{1.35cm}>{\raggedright\arraybackslash}X>{\raggedright\arraybackslash}X}
\toprule
Class & World of events & World of numbers \\
\midrule
GS & everyday happenings; then a bag of counters: certain, maybe, impossible decided by \emph{all}, \emph{some}, \emph{none} & --- \\
CP--CE2 & no chapter in the French-curriculum courses; the reservoir levels used by IB PYP~1--2 reuse the bag method and ask for number-free constructions (``colour the marbles so that picking a red one is certain''), next to chance words on everyday happenings & likelihood line (PYP~2) \\
CM1 & random experiment; outcome; all possible outcomes; event as a set of outcomes, in words (the outcomes that make it happen) & likelihood words on a line; probability from $0$ to $1$; favourable over total; fair games \\
CM2 & two-step experiments: outcomes as ordered results, table and tree & ``the die has no memory''; calculating and comparing in two steps \\
6e & outcomes, events, two-step experiments (review) & fraction, decimal, percentage; repeated trials compared with the probability \\
5e & sample space $U$; event as a subset; singleton; complementary event & $P(E)$; equally likely outcomes; complement rule; relative frequency (titled ``experimental probability'' until 25~October~2026) \\
4e & union, intersection, mutually exclusive events; sample space of a two-step experiment as a set of pairs & probability of an outcome as $P(\{x\})$; sum of the probabilities; fluctuation \\
3e & the language of events (review); counting with two-step experiments & addition relation; simulation; historical problems \\
2nde & sets and events in one table ($\Omega$, $\emptyset$, $\in$, $\cup$, $\cap$, disjoint) & law of probability; conditional probability with two-way tables, weighted trees, inversion; sampling \\
1re & partitions of the sample space; random variable as a function; $(X=x)$ as an event & axioms; conditional probability, total probability, independence; law, expectation \\
Tle & product sample space of $n$ trials; paths with $k$ successes & Bernoulli scheme, binomial law; law of large numbers \\
\bottomrule
\end{tabularx}
\caption{The progression of the French-curriculum courses, read as two parallel strands. GS and the CP--CE2 row describe reservoir levels used by the IB PYP chapters; CP to CE2 have no probability in the French programme or in the French-curriculum courses. From \emph{sixi\`eme} to \emph{seconde}, the frequencies close the chapter (Section~\ref{sec:p5}). Version of 10~October~2026.}
\label{tab:progression}
\end{table}

\paragraph{CM1: the founding order.} The French programme introduces probability in CM1 as ``a first familiarisation with random experiments'', with the objectives of identifying random experiments, identifying all the possible outcomes, using the vocabulary from ``impossible'' to ``certain'', comparing events, and understanding that ``it is not because there are two possible outcomes that each has one chance in two'' \citep{men2025cycle3}. The chapter follows the programme and anticipates by one year two objectives that the programme places in CM2: identifying ``all the outcomes realising an event'' and expressing a probability as ``$a$ chances in $b$'' \citep{men2025cycle3}. It adds one thing that neither year names: the definition of an event as a set of outcomes, with the example ``find the outcomes that correspond to rolling an even number''.

\paragraph{CM2 and \emph{sixi\`eme}: products and forms.} CM2 is the year of the first non-trivial sample space, the ordered pairs of a two-step experiment, and of independence as ``the die has no memory''; the programme asks for both \citep{men2025cycle3}. \emph{Sixi\`eme} reviews the vocabulary and devotes most of its new content to the three written forms of a probability---fraction, decimal, percentage---and to the comparison of a series of trials with the computed probability. The programme makes the passage from ``$a$ chances in $b$'' to the number $\frac{a}{b}$ ``a major objective'' of the year, and states that pupils are not expected to use the specific vocabulary of probability on their own \citep{men2025cycle3}; the chapter keeps that vocabulary in use.

\paragraph{\emph{Cinqui\`eme} and \emph{quatri\`eme}: the algebra named.} The \emph{cinqui\`eme} chapter is the first to call its first part \emph{The Algebra of Events}, with the sample space $U=\{x_1,\dots,x_n\}$, events as subsets, singletons and complements; the programme of that year asks only for the vocabulary ``random experiment, outcome, event'' \citep{men2026cycle4}, so that the chapter anticipates by one year the set notation that the programme places in \emph{quatri\`eme}. \emph{Quatri\`eme} adds union, intersection and mutually exclusive events, and makes the probability of an outcome that of its singleton.

\paragraph{\emph{Troisi\`eme} to \emph{premi\`ere}: the rules and the axioms.} \emph{Troisi\`eme} states the addition relation in the programme's form $P(A\cup B)+P(A\cap B)=P(A)+P(B)$ and introduces simulation; \emph{seconde} gathers the vocabulary of sets and events in a single table and introduces conditional proportions; \emph{premi\`ere} gives the axioms, partitions, conditional probability and random variables, with the explicit statement of the principle quoted in Section~\ref{sec:p1}.

\paragraph{\emph{Terminale}: products and laws.} The \emph{terminale} chapters build the sample space of $n$ trials as a Cartesian product, model independence as a choice, and define the binomial law by counting the outcomes of the event ``$k$ successes''.

\section{What the analysis reveals about the chapters}\label{sec:audit}

The lessons state the principles; the exercises show what pupils are asked to do. Table~\ref{tab:exercises} codes the exercises of each French-curriculum chapter by the competence of their root block, with the categories of Section~\ref{sec:intro}; Table~\ref{tab:written} examines the written exercises, which are corrected against a rubric.

\begin{table}[ht]
\centering\small
\begin{tabular}{lrrrrrrrr}
\toprule
Class & events & modelling & likelihood & numerical & problems & total & events & events \\
 & & & words & & & & share & + modelling \\
\midrule
CM1  & 21 & 12 & 17 &  16 &  6 &  72 & 29\% & 46\% \\
CM2  & 16 &  8 &  4 &  10 &  4 &  42 & 38\% & 57\% \\
6e   & 13 &  8 &  7 &  32 & 11 &  71 & 18\% & 30\% \\
5e   & 24 &  0 & 19 &  46 & 10 &  99 & 24\% & 24\% \\
4e   & 32 &  0 &  0 &  36 &  5 &  73 & 44\% & 44\% \\
3e   & 32 &  0 &  0 &  60 &  7 &  99 & 32\% & 32\% \\
2nde & 36 &  0 &  0 & 105 &  7 & 148 & 24\% & 24\% \\
\bottomrule
\end{tabular}
\caption{Exercises of the French-curriculum probability chapters, coded by the competence of their root block in the reservoir (the same block receives the same code in every class). Categories as in Section~\ref{sec:intro}. English half of the sources, version of 10~October~2026.}
\label{tab:exercises}
\end{table}

\begin{table}[ht]
\centering\small
\begin{tabular}{>{\raggedright\arraybackslash}p{6.4cm}rrrrrrr}
\toprule
 & CM1 & CM2 & 6e & 5e & 4e & 3e & 2nde \\
\midrule
written exercises & 10 & 5 & 14 & 19 & 13 & 22 & 31 \\
\quad of which: rubric assesses the construction of the event & 1 & 0 & 0 & 5 & 4 & 5 & 1 \\
exercises whose statement asks for an event or a sample space as a set & 0 & 0 & 0 & 8 & 5 & 0 & 3 \\
\bottomrule
\end{tabular}
\caption{Written exercises and the event. A rubric is counted when it rewards or penalises the construction of the event---for instance ``Penalise a missing outcome or a list without the set braces'' (\emph{cinqui\`eme}) or ``Penalise: an answer [\ldots] with no event listed first'' (\emph{quatri\`eme}). Statements are counted when they ask to ``write'' or ``determine'' an event or a sample space ``as a set''.}
\label{tab:written}
\end{table}

The confrontation brings out the following points.
\begin{enumerate}[label=(\alph*)]
  \item\label{it:balance} \emph{The balance of the exercises.} Every chapter has exercises on events, and the \emph{quatri\`eme} chapter gives each operation---union, intersection, mutually exclusive events, complement---its own set. The share of event exercises varies from $18\%$ (\emph{sixi\`eme}) to $44\%$ (\emph{quatri\`eme}); with the modelling exercises of the primary classes (equally likely outcomes, fair games, independence), the share reaches $46\%$ in CM1 and $57\%$ in CM2. In \emph{sixi\`eme}, the subsection on the three forms of a probability (fraction, decimal, percentage) has fourteen exercises; part of them start from a draw and count the outcomes of an event before converting, the others convert a given number, in line with the programme's ``major objective'' for the year \citep{men2025cycle3}. Exercise counts do not measure time. They show that in CM1 and CM2 events and modelling make up about half of the exercises ($46\%$ and $57\%$), whereas from \emph{sixi\`eme} on, \emph{quatri\`eme} excepted, more than half ask for likelihood words or numbers ($55\%$ to $71\%$).
  \item\label{it:answers} \emph{The event in the calculation exercises.} In primary school the design keeps events first without set notation, and the exercises follow it: event exercises ask pupils to identify the outcomes that make an event happen (in CM1, for instance, ``Identifying outcomes for die-rolling events''), and the model answers of the calculation exercises name the event and count its outcomes (``The number of outcomes for the event `choosing a red ball' is 2''). From \emph{cinqui\`eme}, set notation enters the written exercises (Table~\ref{tab:written}): statements ask pupils to ``write $E$ and $F$ as sets of outcomes'' or to ``determine the event $E$ [\ldots] as a set of pairs, then compute $P(E)$'', and rubrics penalise ``describing an event without writing it as a set'' or ``a computation [\ldots] with no event described''. What the audit questions is the format of the automatically corrected calculation exercises: at every level, their answer field follows $P(\text{event})=$ and records the number alone, so that the construction of the event, which the method prescribes as the first step, is shown in the model answer but not asked of the pupil. An answer in two steps---first the outcomes of the event, as a list in primary school and as a set from \emph{cinqui\`eme}, then its probability---would assess the order that the design prescribes.
  \item\label{it:everyday} \emph{Events without a stated experiment.} The likelihood line shared by CM1 to \emph{cinqui\`eme} has five examples; two are events of experiments (``rolling a die and getting a 3'', ``tossing a coin and getting heads''), three are everyday happenings (``riding a dinosaur'', ``drinking water at school today'', ``the sun will rise tomorrow''), and the GS chapter opens in the same way. The CM1 programme itself uses only events of experiments \citep{men2025cycle3}. Everyday happenings are legitimate objects of likelihood judgements, and each could be modelled as an event of a two-outcome experiment; but the chapters do not state such a model, so that the definition of an event given a few lines earlier is not applied to them. A design faithful to Principle~\ref{pr:two} would draw the line's examples from an experiment in view, or state the model.
  \item\label{it:nomodel} \emph{Numbers without a universe.} Several examples give a probability without naming any sample space: ``Farid has a $0.8$ ($80\%$) chance of finishing his homework on time tonight'' (\emph{cinqui\`eme} to \emph{seconde}), ``a $40\%$ chance of rain tomorrow'' (\emph{quatri\`eme}, \emph{troisi\`eme}), ``the probability that a randomly selected student participates in singing is $0.4$'' (\emph{troisi\`eme}, \emph{seconde}). The last has a universe---the students of the school---that is not named. Such examples are useful because the measure is not obtained by counting; they would serve the principle better if the universe were stated even when the measure is given.
  \item\label{it:rules} \emph{Derivations shown in the \emph{coll\`ege}, hidden in the lyc\'ee.} The complement rule is proved from the partition of the sample space in \emph{cinqui\`eme}, \emph{quatri\`eme} and \emph{troisi\`eme}, and the addition relation from the double counting of the intersection in \emph{troisi\`eme}; in \emph{seconde} both proofs are hidden, and in \emph{premi\`ere} the proofs of the complement rule and of the sum over a partition are hidden too, while monotonicity is proved (Section~\ref{sec:p4}). The pattern is the reverse of what one might expect: the rules are carried over from events for younger pupils and stated as facts for older ones.
  \item\label{it:conditional} \emph{One view of conditioning.} Conditional probability is presented as a proportion in a new population, which is correct and suited to counts; the second view, a new probability on the same events, is not stated (Section~\ref{sec:p5}).
  \item\label{it:frequency} \emph{A frequency called a probability.} The block that introduces frequencies from \emph{cinqui\`eme} to \emph{seconde} is titled \emph{Experimental Probability (Relative Frequency)}. Its definition calls the number ``an estimate'', and its opening problem writes $P(\text{base down})\approx\frac{15}{50}$, which keeps the letter $P$ for the model; but the name suggests a second kind of probability rather than a statistic that estimates one. The \emph{sixi\`eme} lesson avoids it and speaks of ``the part of the trials''. Naming the object \emph{frequency} throughout, as the \emph{seconde} chapter does when it writes $f$ for the frequency observed in a sample, would extend Principle~\ref{pr:two} to a third object: the frequency of an event in a series of trials is neither the event nor its probability.
  \item\label{it:dalembert} \emph{A misattribution.} The \emph{troisi\`eme} historical problem attributed to d'Alembert counts ``three possible results---$0$, $1$ or $2$ heads''; this is a later, folklore version of the error. D'Alembert's own argument, in the \emph{Encyclop\'edie}, concerns the three cases of Figure~\ref{fig:dalembert} \citep{dalembert1754,henry2001croix}.
  \item\label{it:early} \emph{The early years.} The French programme of cycle~2 contains no probability \citep{men2024cycle2}, and the French-curriculum courses have none before CM1. The reservoir levels that the IB chapters use for ages five to seven contain the bag method of Figure~\ref{fig:bag} and number-free constructions of events, next to chance words applied to everyday happenings (``it is impossible for a dog to meow like a cat''), and some likelihood items that treat necessity as chance (``how likely is it that a circle is round? Certain''). The germ of the event-first order is present at these ages, but not organised; Section~\ref{sec:cycle2} proposes how it could be.
\end{enumerate}
Points \ref{it:answers} to \ref{it:dalembert} can be corrected within the existing blocks: a calculation answered in two steps, the outcomes of the event before the number; an example that names its universe; a proof that is displayed rather than hidden; a sentence on $P_A$; a frequency named as such; a corrected attribution. Point \ref{it:balance} suggests more exercises on events, which the reservoir partly contains, many of them near-copies of exercises already linked at another level. Point \ref{it:early} is a design decision.

\paragraph{After the analysis.} Several issues identified by the audit were addressed, wholly or partly, on 16~October~2026, in the reservoir and therefore in the French-curriculum chapters that link its blocks. Eight automatic exercises in CM1, twelve in CM2 and twelve in \emph{sixi\`eme} now ask pupils to enter the number of outcomes in the event and the total number of possible outcomes before the probability (point~\ref{it:answers}); these fields assess counts rather than an explicit list of the event's outcomes. From \emph{cinqui\`eme} to \emph{seconde}, the automatic calculation exercises still record the number alone. The likelihood line now says that its everyday happenings are not events of an experiment whose outcomes are listed (point~\ref{it:everyday}). The examples of Farid's homework and of tomorrow's rain now state a model with two outcomes, such as $U=\{\text{on time},\text{late}\}$ (point~\ref{it:nomodel}). The second view of conditioning is stated in the \emph{premi\`ere} chapter as a proposition, \emph{Conditioning Changes the Probability, Not the Events}, whose proof is hidden there (point~\ref{it:conditional}). The \emph{troisi\`eme} exercise on d'Alembert is rewritten on the argument of ``Croix ou pile'' (point~\ref{it:dalembert}). The proofs of the \emph{seconde} remain hidden (point~\ref{it:rules}). The IB chapters were revised at the same time (Section~\ref{sec:ib}). On 25~October~2026, after a reading of the \emph{troisi\`eme} chapter that an external model reviewed, the block on frequencies was renamed \emph{Relative Frequency} (\emph{La fr\'equence}), ``experimental probability'' being kept as a second name (point~\ref{it:frequency}); the opening text on equally likely outcomes now presents equiprobability as a hypothesis of the model, made reasonable by symmetry, as the programme requires; and every chapter from CM2 to \emph{premi\`ere} received a common skeleton: the language of events, then the probability, then, up to \emph{troisi\`eme}, a last section \emph{From Chance to Experiment} with a single subsection (the \emph{seconde} and \emph{premi\`ere} keep their last sections, \emph{From the Model to the Sample} and \emph{Estimating by Simulation}). In \emph{quatri\`eme} and \emph{troisi\`eme} that last section opens on a new example comparing the diagram of the frequencies observed in sixty throws of a die with the diagram of the model's probabilities. In CM2, ``the die has no memory'' became a remark in the last section, and the method for two-step experiments now justifies the equally likely cells by fair steps carried out under the same conditions, no longer by the die's lack of memory; the outcome of a two-step experiment is called a pair, as in the video that introduces it. The model answers on bags and spinners from CE2 to CM2 now number the balls, or the equal parts of the spinner, say that each has the same chance, and only then count: ``the event `choosing a red ball' is the set of balls 1 and 2''; a figure puts the balls of the event in a dashed hoop. The set is written in words, braces still arriving in \emph{cinqui\`eme}, and the CE2 lesson used in CM1 now states once how a bag is drawn from: marbles alike, mixed, drawn without looking, each marble being an outcome. The IB chapters received the same changes on 25~October~2026: in PYP~4 the language of events now precedes the numbers; in PYP~5 the outcome of a two-step experiment is the pair, and the die without memory is a closing remark; in MYP~3 the table of outcomes comes before the probabilities; and relative frequency moves to the last section of the MYP~2, MYP~4 and DP~1 chapters.

\paragraph{A comparison with other curricula.} The figures above have no absolute meaning, but the place of events in other curricula gives a point of comparison. In the NCTM \emph{Principles and Standards}, the probability of an event is represented by a number from $0$ to $1$ in grades 3--5, while ``the concepts of sample space'' are expected only in grades 9--12 \citep{nctm2000}. The Common Core defines probability in grade~7 as ``a number between 0 and 1 that expresses the likelihood of the event occurring'', asks pupils to identify ``the outcomes in the sample space which compose the event'' for compound events, and asks them to describe ``events as subsets of a sample space'' only in high school \citep{ccss2010}. The GAISE~II framework does not use the phrase ``sample space'' \citep{bargagliotti2020}. The French programme places events as sets, with their operations, in \emph{quatri\`eme} \citep{men2026cycle4}; the IB Diploma Programme names ``sample space ($U$) and event'' on the same line as the formula $P(A)=n(A)/n(U)$ \citep{ibo2019aa}. Against these documents, the corpus is unusual in defining an event as a set of outcomes from CM1, at age nine; its exercises, by contrast, are closer to the common pattern than its lessons. A comparison at the level of textbooks, with the same coding, remains to be done.

\section{The early years: events before numbers in CP and CE1}\label{sec:cycle2}

The French programme leaves chance out of cycle~2, and the IB PYP, which includes it from the first phase, treats it mainly through likelihood vocabulary applied to ``daily events'' \citep{ibo2009pyp}. The thesis of this paper suggests a different entry for ages six to eight: not numbers, which the programme rightly postpones, and not likelihood words alone, but the \emph{construction of events without numbers}. The kindergarten method of Figure~\ref{fig:bag} and the PYP~2 task ``colour the marbles so that picking a red one is certain'' show that the core of it is within reach of young children. Four kinds of task organise it, each using only objects that the pupil can see and move.
\begin{enumerate}[label=(\arabic*)]
  \item \emph{What can come out?} Before a draw from a bag, a roll of a die or a spin, the pupil lists or shows everything that can come out---the counters, the faces, the sectors. This is the sample space, built by hand.
  \item \emph{Which ones make it happen?} For a sentence about the result (``the counter is red'', ``the die shows more than four''), the pupil puts in a hoop the outcomes that make the sentence true. The hoop is the event; the sentence is its name.
  \item \emph{All, some, none.} The pupil decides whether the event is certain, possible or impossible by comparing the hoop with the whole: all the outcomes are in the hoop, some are, none is.
  \item \emph{Not, or, and.} With two hoops, the pupil finds the outcomes that are \emph{not} in a hoop, those in one hoop \emph{or} the other, those in both; the words come before any symbol, and the pupil checks, for instance, that ``red or not red'' is certain and ``red and not red'' impossible.
\end{enumerate}
None of these tasks requires a number, a fraction or a formula. The PYP already asks learners, in its second phase, to ``represent the relationship between objects in sets using tree, Venn and Carroll diagrams'' \citep{ibo2009pyp}; the French programme of cycle~2 sorts objects by a criterion---solids, for instance---but does not use Venn or Carroll diagrams \citep{men2024cycle2}. What the tasks add in both cases is the reference to a random experiment: the objects sorted are possible results, and the regions are events.

\paragraph{An a priori analysis of the two-hoop task.} The richest of the four is task~(4). Its didactic variables, in the sense of the theory of didactical situations \citep{brousseau2001}, are: the number of outcomes (from four to eight counters); their distinguishability (counters that differ by one attribute, colour, or by two, colour and shape, so that two sentences can cut across each other); the way the events are given (by a sentence about an attribute, by a list, or by a hoop already filled); the overlap of the two events (disjoint, nested, crossing); and the presence or absence of a symmetry that makes the outcomes equally likely. The expected procedures are: placing each counter by checking the sentence (a correct procedure that builds the event by extension); placing only the ``typical'' counters (an incomplete event); and, for ``or'', placing only the counters that satisfy exactly one sentence (the exclusive reading). The anticipated errors are those of \citet{way2003}: judging ``impossible'' when the counters of a colour are few rather than absent, and ``certain'' when they are many rather than all. The feedback comes from the situation, with a limit that is itself instructive. One draw outside the hoop refutes ``certain'', and one draw inside it refutes ``impossible''; but no number of confirming draws establishes either judgement, which can only be settled by inspecting all the outcomes---the passage from the draw to the set that the task is meant to provoke. A judgement of ``possible'' is confirmed by a single draw inside the hoop and refuted by none. With crossing events and two attributes, the pupil can see that the hoop ``red or round'' holds fewer counters than the red ones and the round ones counted separately, because the red round counters would be counted twice---the double counting of Figure~\ref{fig:transport}, before any number.

Two cautions apply. The first is the one of Section~\ref{sec:background}: these tasks are not set theory for its own sake; the hoops are tools for describing a random experiment, the words \emph{not}, \emph{or}, \emph{and} are ordinary words used with a precise meaning, and what is assessed is the judgement about chance. The second is that the proposal is a design hypothesis: at the time of the analysis, the reservoir contained the bag method and the construction of a certain event, but not tasks (2) and (4) for these ages. The revision of 16~October~2026 added hoop tasks and work on \emph{not}, \emph{or}, \emph{and} to the IB PYP~0 to PYP~2 chapters (Section~\ref{sec:ib}); their effect on learning remains a hypothesis, and Hypothesis~\ref{h:early} states how it could be measured.

\section{Developing the separation in the IB curriculum}\label{sec:ib}

\paragraph{The official documents.} The IB documents place the two worlds at different stages (Table~\ref{tab:ib}). The PYP scope and sequence describes chance mainly on the side of likelihood: learners ``discuss chance in daily events (impossible, maybe, certain)'', then describe likelihood with words; in a later phase they learn that probability is based on experimental events and can be expressed numerically, use tree diagrams to express probabilities ``using simple fractions'', and ``use probability to determine mathematically fair and unfair games and to explain possible outcomes''; finally they express probability on a scale from $0$ to $1$ and predict it theoretically \citep{ibo2009pyp}. The same phase that introduces chance words asks learners to represent ``the relationship between objects in sets using tree, Venn and Carroll diagrams'', without connecting the two. The document stresses that its phases ``are not always linear or age related'', so that any mapping of phases onto PYP years, including the one of Table~\ref{tab:ib}, is an approximation. The MYP framework lists ``sample spaces'' and ``probability of simple events'' in MYP~1--3, and ``sets, including notation and operations up to three sets'' with ``probability with Venn diagrams, tree diagrams and sample spaces'', mutually exclusive and combined events in MYP~4--5 \citep{ibo2020myp}. The Diploma Programme guide lists the notation of sets and their union and intersection as prior learning, and states the concepts of ``trial, outcome, equally likely outcomes, relative frequency, sample space ($U$) and event'' on the same line as $P(A)=n(A)/n(U)$; it recommends Venn diagrams, tree diagrams, sample space diagrams and tables of outcomes, notes ``the non-exclusivity of `or'\,'', and adds that problems ``can be solved [\ldots] without explicit use of formulae'' \citep{ibo2019aa}. It defines mutually exclusive events by ``$P(A\cap B)=0$'', a numerical criterion that conflates the impossible event with an event of probability $0$---the distinction drawn in Section~\ref{sec:maths}; the set criterion $A\cap B=\emptyset$ belongs to the world of events. These guides are being replaced: the IB announces a new mathematics curriculum taught from 2027 \citep{ibo2026dp}, so that the proposals below are addressed to its implementation as much as to the 2019 guides.

\begin{table}[ht]
\centering\footnotesize
\begin{tabularx}{\textwidth}{>{\raggedright\arraybackslash}p{1.5cm}>{\raggedright\arraybackslash}X>{\raggedright\arraybackslash}X}
\toprule
Stage & Official IB documents (paraphrased) & CommeUnJeu IB chapters (after the revision) \\
\midrule
PYP, early phases & chance in daily events; likelihood words; (data handling: objects sorted in sets with tree, Venn and Carroll diagrams) & PYP~0: the bag before the words (all, some, none); the red counters put in a hoop; PYP~1: what can come out of a bag, a bag coloured for a given event; PYP~2: \emph{not}, \emph{or}, \emph{and} on a bag, then the likelihood line \\
PYP, later phases & probability from experiments, expressed numerically; tree diagrams with simple fractions; fair games and possible outcomes; scale $0$--$1$; theoretical prediction & PYP~3--5: outcomes and events as sets of outcomes, then words and numbers; two-step experiments; fair games; experiments at the end of the chapter \\
MYP 1--3 & sample spaces; probability of simple events; theoretical and experimental probability & MYP~1: \emph{Outcomes and Events}, MYP~2: \emph{Sample Spaces and Events}, both \emph{before} likelihood words and numbers; MYP~3: events as sets first, and a separate chapter on sets and Venn diagrams \\
MYP 4--5 & sets and operations; probability with Venn diagrams, trees and sample spaces; mutually exclusive and combined events; (extended) conditional probability & MYP~4: events in the language of sets, then their probabilities; MYP~5: Venn diagrams with three sets, the algebra of events, the axioms, the rules carried over from the events \\
DP & sets as prior learning; sample space $U$ and event; $P(A)=n(A)/n(U)$; Venn, tree and sample-space diagrams; mutually exclusive as $P(A\cap B)=0$; conditional and independent events & a chapter on sets and Venn diagrams; a probability chapter that opens on the world of events before the numbers; conditional probability in a chapter of its own \\
\bottomrule
\end{tabularx}
\caption{The two worlds in the IB curriculum: what the official documents prescribe \citep{ibo2009pyp,ibo2020myp,ibo2019aa}, in paraphrase, and what the IB chapters of the platform contain after their revisions of 16 and 25~October~2026 (for PYP~4, the move of the experiments to the end of the chapter dates from 22~October). The mapping of PYP phases onto years is approximate.}
\label{tab:ib}
\end{table}

\paragraph{The chapters of the platform.} The IB chapters of the platform are assembled from the same reservoir as the French ones. In the version of 10~October~2026, they inherited its order from PYP~3 on: outcomes, then events as sets of outcomes, then words, then numbers. They departed from it at the start of the MYP: the MYP~1 chapter, which shares most of its blocks with the French \emph{sixi\`eme}, opened with \emph{Describing Probability} (likelihood words, then probability as a number from $0$ to $1$) and only then turned to \emph{Outcomes and Events}, the reverse of the French order; the MYP~2 chapter also opened with likelihood words and the probability scale, as review, before the sample space. From MYP~3 on, a separate chapter on sets and Venn diagrams, present again in DP~1, gives the algebra of events its own home. The chapters were revised on 16~October~2026, after this analysis, along the proposals below, and again on 25~October, when the changes made to the French chapters were carried over; the right column of Table~\ref{tab:ib} describes them as revised.

\paragraph{Proposals.} Four changes would carry the separation through the IB curriculum, each within the official documents. (1) In the early PYP, connect the two strands that the scope and sequence keeps apart---sorting objects in Venn and Carroll diagrams, and the language of chance---with the tasks of Section~\ref{sec:cycle2}, which turn ``certain'' and ``impossible'' into full and empty regions and prepare the ``possible outcomes'' of the later phases and the ``sample space'' of MYP~1. (2) In MYP~1 and MYP~2, restore the order of the French \emph{sixi\`eme} and \emph{cinqui\`eme}: outcomes, sample space and events before likelihood and numbers. (3) In MYP~3--4, make the chapter on sets the declared prerequisite of probability and use its vocabulary in the probability chapter, as the dictionary table of the French \emph{seconde} does. (4) In the DP, define mutually exclusive events by $A\cap B=\emptyset$ rather than by $P(A\cap B)=0$, and teach the items that belong to the world of events---sample space diagrams, tables of outcomes, the inclusive ``or''---as such, before the formulas, which is what the guide's remark on solving problems ``without explicit use of formulae'' invites.

\paragraph{The proposals in the platform's chapters.} The revision of 16~October~2026 applies the four proposals to the IB chapters of the platform. (1) The PYP~0 chapter shows the bag before the judgement \emph{certain}, \emph{maybe} or \emph{impossible}, and has pupils put the red counters in a hoop; PYP~1 opens with \emph{What Can Come Out?} and asks pupils to colour the marbles of a bag so that an event is certain, possible or impossible; PYP~2 opens with \emph{Not, Or, And}, before the likelihood words. (2) MYP~1 opens with \emph{Outcomes and Events}, and MYP~2 with \emph{Sample Spaces and Events}; both place them before likelihood words and numerical probabilities. (3) The probability chapters of MYP~3 and MYP~4 open with events written as sets, and, after a section on Venn diagrams with three sets, MYP~5 follows the proposed sequence: \emph{The Algebra of Events}, \emph{The Axioms of Probability}, \emph{The Rules, Carried Over from the Events}. (4) In DP~1, where the chapters already defined mutually exclusive events by $E\cap F=\emptyset$, the probability chapter opens with \emph{The World of Events}---sample spaces, \emph{not}, \emph{or}, \emph{and}, partitions---before the numbers, and conditional probability has a chapter of its own. Whether these changes improve learning is, like the rest of the design, a hypothesis.

\section{Anticipated misconceptions}\label{sec:misconceptions}

Table~\ref{tab:misconceptions} lists the misconceptions that the design addresses, the level at which each arises---the events, the weights given to outcomes, or the numbers---and the response of the progression.

\begin{table}[ht]
\centering\footnotesize
\begin{tabularx}{\textwidth}{>{\raggedright\arraybackslash}p{3.7cm}>{\raggedright\arraybackslash}p{1.3cm}>{\raggedright\arraybackslash}X>{\raggedright\arraybackslash}p{1.3cm}}
\toprule
Misconception & Level & Response of the progression & Class \\
\midrule
The event is its probability (``the event is one half'') & both & two names; $P$ as a function from events to numbers & 5e, 1re \\
A question about chance read as a prediction of the next outcome \citep{konold1989} & events & list the outcomes, then name the event, before any likelihood & CM1 \\
Colours counted instead of balls; unequal sectors counted as equal (equiprobability bias, \citealp{lecoutre1992}) & weights & ``count the balls, not the colours''; two spinners compared & CM1 \\
``A five and a six'' as likely as ``two sixes''; d'Alembert's three cases; totals of dice counted as equally likely \citep{fischbein1997} & weights & refinement into equally likely ordered pairs, by table or tree; historical problems & CM2, 4e, 3e \\
$P(A\cup B)=P(A)+P(B)$ for overlapping events; probabilities above $1$ & both & the intersection counted twice, with a proof & 3e \\
Mutually exclusive taken for complementary & events & a remark with two disjoint events that do not cover $\Omega$ & 4e \\
$P(A\cap B)>P(A)$ (conjunction fallacy, \citealp{tversky1983}) & both & monotonicity $A\subseteq B\Rightarrow P(A)\le P(B)$, with a proof; the check ``$P(A\cap B)$ can exceed neither $P(A)$ nor $P(B)$'' & 3e, 2nde, 1re \\
A result ``due'' after a run (gambler's fallacy) & weights & ``the die has no memory'': independence of repetitions & CM2 \\
$P_A(B)$ confused with $P_B(A)$ or with $P(A\cap B)$ & both & the known event as a new population; screening tests & 2nde, 1re \\
An everyday happening treated as an event without a model & events & design observation (Section~\ref{sec:audit}) & --- \\
\bottomrule
\end{tabularx}
\caption{Anticipated misconceptions, the level at which they arise, and the response of the progression. Version of 10~October~2026.}
\label{tab:misconceptions}
\end{table}

\section{Hypotheses for an empirical evaluation}\label{sec:hypotheses}

No data are reported in this paper; its claims about learning are hypotheses. A study would compare pupils taught with the progression and pupils taught with a curriculum that introduces probabilities through likelihood and ``favourable over total'' without an explicit construction of events, matched for time and content, and would record the sequence actually taught in each group. Three hypotheses bear on the central claim.

\begin{hypothesis}[The object]\label{h:object}
Pupils taught events before probabilities more often describe an event as a set of outcomes, or by a sentence together with its outcomes, and less often as a number. Measures, at the end of \emph{cinqui\`eme} and \emph{troisi\`eme}: ``What is an event? Give an example for one roll of a die''; judging statements such as ``$A=\frac12$'' and ``$P(A)=\{2,4,6\}$''; coded as set or description of outcomes, number, other.
\end{hypothesis}

\begin{hypothesis}[The universe]\label{h:universe}
The same pupils make fewer errors in which outcomes are wrongly weighted. Measures, at the end of CM2 and \emph{quatri\`eme}: a bag with unequal numbers of balls of each colour; a spinner with unequal sectors; ``a five and a six'' against ``two sixes''; d'Alembert's two tosses; coded for whether the answer uses a correctly weighted model---equally likely refined outcomes, or coarse outcomes with correct weights---and for the equiprobability bias.
\end{hypothesis}

\begin{hypothesis}[The algebra and the rules]\label{h:algebra}
The same pupils solve more combined-event problems correctly, give fewer probabilities above $1$, and give fewer answers with $P(A\cap B)>P(A)$. Measures, at the end of \emph{troisi\`eme}: union of overlapping events, complement, conjunction items; coded for the addition error, the conjunction fallacy and the use of a diagram of outcomes.
\end{hypothesis}

Four secondary hypotheses concern the extensions and the place of frequencies.

\begin{hypothesis}[Early years]\label{h:early}
Pupils who have done the tasks of Section~\ref{sec:cycle2} in CP and CE1 identify the outcomes of an event, and judge certain and impossible events, more accurately in CM1 than pupils who have learnt chance words only.
\end{hypothesis}

\begin{hypothesis}[Conditional probability]\label{h:cond}
Pupils taught both views of conditioning, including $P_A$ as a new probability on the same events, confuse $P_A(B)$ with $P(A\cap B)$ less often.
\end{hypothesis}

\begin{hypothesis}[Frequencies]\label{h:frequency}
Pupils who meet frequencies at the end of the chapter, as estimates of the probability of an event already built, confuse a frequency with a probability less often: fewer conclude from a short series that a die is loaded, or expect a series to reproduce the probability exactly. Measures, at the end of \emph{cinqui\`eme} and \emph{troisi\`eme}: the twelve rolls without a six; a cone thrown $50$ times; coded for the distinction between the observed frequency and the probability of the model.
\end{hypothesis}

\begin{hypothesis}[Random variables]\label{h:rv}
Pupils who define $(X=k)$ as an event compute laws from the outcomes more often than from memorised formulas.
\end{hypothesis}

Each hypothesis should be tested with a non-inferiority margin on ordinary computational items: the claim is that the event-first order improves understanding without loss on calculation.

\section{Limitations}\label{sec:limitations}

The progression and the analysis are the work of one author, who has a commercial interest in the materials; the audit is a self-audit. The counts were produced by scripts over the source mark-up, with the help of an AI assistant, and checked by reading; an independent referee's check of a first version found errors in the inventory---proofs attached to propositions had been missed and the coding was inconsistent across classes---which the present version corrects by coding each exercise by its root block and by examining proofs and rubrics. The coding sheet is available with the source documents, and a second coder would strengthen the counts. The corpus is a dated version of chapters that change; the counts describe what the chapters ask, not what teachers assign or pupils do, and no data on classroom practice are available. The comparison with other curricula is made at the level of programmes, not of textbooks. The IB analysis rests on the official documents, which are paraphrased and are being replaced, and on the platform's chapters; it does not describe IB classrooms. The thesis concerns finite probability; the continuous case of the \emph{terminale}, in which the events met are intervals and their unions and probabilities are computed as integrals of a density, extends the separation but is not analysed. The everyday, subjective use of probability for events that are not outcomes of a stated experiment is legitimate; the paper argues only that it should not be confused with the events of a model. Finally, the objection of the ``modern mathematics'' reform is answered by argument, not by evidence: whether early work on events as sets avoids the failures of that reform is itself an empirical question.

\section{Conclusion}\label{sec:conclusion}

An event is a set of outcomes; its probability is a number. The paper has argued that this distinction is not a matter of vocabulary but the structure of the subject; that several famous errors in the history of probability concern the weighting of outcomes, and that their historical remedy was a refinement of events into equally likely outcomes; and that research on learners locates many of their difficulties at the same place. It has proposed five principles that order the teaching accordingly---two objects, the event first, the modal before the numerical, the algebra and its rules, the universe as a model, with the frequency at the end of the chapter as a statistical estimate---and shown that the French-curriculum courses of the platform realise them in their lessons from CM1 to the \emph{terminale}, with texts that state the principle in so many words and, in the \emph{coll\`ege}, proofs that carry the rules over from events. In the version analysed, the audit has also shown where the realisation stops short: exercises in which events are a minority outside \emph{quatri\`eme}, calculation exercises whose automatic answer is the number alone, examples without a universe, proofs hidden in the lyc\'ee, and early years in which the event-first order exists only in germ. Several of these points were addressed, wholly or partly, after the analysis (Sections~\ref{sec:audit} and~\ref{sec:ib}). Two extensions follow: number-free tasks for CP and CE1 in which pupils build events, with an a priori analysis of the richest of them; and changes to the IB curriculum that connect the PYP's sorting with its chance, restore the event-first order in MYP~1 and MYP~2, and define mutually exclusive events by their intersection; the platform's IB chapters were revised along these lines after the analysis. Whether the order improves learning is an empirical question; the paper states it as hypotheses. Its claim is narrower and, we hope, firm: before learning to compute the probability of an event, pupils should learn to build, identify and handle the event, because that is what a probability model is made of.

\section*{Competing interests}
The author is the founder of CommeUnJeu, which publishes the materials analysed in this paper and offers paid services on its website.

\section*{Acknowledgements}
The author used AI language tools while preparing this paper: Claude (Anthropic) assisted with drafting, with the inventory and counting of the corpus and with the checking of bibliographic references, and acted, in a separate instance, as an independent referee; ChatGPT (OpenAI) was used as a critical reader. The author designed the progression, checked every reference against its source, and takes full responsibility for the content.

\bibliographystyle{abbrvnat}
\bibliography{probability_refs}

\appendix
\section{Public resources}\label{app:resources}

All addresses begin with \texttt{https://commeunjeu.org/cours/}. In the French slugs, \texttt{grade-}$N$ denotes the $N$-th year of schooling counted from CP.

\begin{table}[ht]
\centering\small
\begin{tabularx}{\textwidth}{@{}l>{\raggedright\arraybackslash}X@{}}
\toprule
Class & French-curriculum chapter \\
\midrule
CM1 & \nolinkurl{fr-grade-4-probability/} \\
CM2 & \nolinkurl{fr-grade-5-probability/} \\
6e & \nolinkurl{fr-grade-6-probability/} \\
5e & \nolinkurl{fr-grade-7-probability/} \\
4e & \nolinkurl{fr-grade-8-probability/} \\
3e & \nolinkurl{fr-grade-9-probability/} \\
2nde & \nolinkurl{fr-grade-10-probability-and-sampling/} \\
1re & \nolinkurl{fr-grade-11-speciality-conditional-probability/} \\
1re & \nolinkurl{fr-grade-11-speciality-random-variables/} \\
Tle & \nolinkurl{fr-grade-12-speciality-combinatorics-and-counting/} \\
Tle & \nolinkurl{fr-grade-12-speciality-independent-trials-bernoulli-scheme-and-binomial-distribution/} \\
Tle & \nolinkurl{fr-grade-12-speciality-sums-of-random-variables-concentration-and-the-law-of-large-numbers/} \\
\bottomrule
\end{tabularx}
\caption{Probability chapters of the French-curriculum courses (October 2026).}
\label{tab:addresses-fr}
\end{table}

\begin{table}[ht]
\centering\small
\begin{tabularx}{\textwidth}{@{}l>{\raggedright\arraybackslash}X@{}}
\toprule
IB year & IB chapter \\
\midrule
PYP 0 & \nolinkurl{ib-pyp-0-sorting-picture-graphs-and-chance/} \\
PYP 1 & \nolinkurl{ib-pyp-1-chance/} \\
PYP 2 & \nolinkurl{ib-pyp-2-chance/} \\
PYP 3 & \nolinkurl{ib-pyp-3-chance/} \\
PYP 4 & \nolinkurl{ib-pyp-4-probability/} \\
PYP 5 & \nolinkurl{ib-pyp-5-probability/} \\
MYP 1 & \nolinkurl{ib-myp-1-probability/} \\
MYP 2 & \nolinkurl{ib-myp-2-probability/} \\
MYP 3 & \nolinkurl{ib-myp-3-sets-and-venn-diagrams/} \\
MYP 3 & \nolinkurl{ib-myp-3-probability/} \\
MYP 4 & \nolinkurl{ib-myp-4-probability/} \\
MYP 5 & \nolinkurl{ib-myp-5-probability/} \\
MYP 5 & \nolinkurl{ib-myp-5-conditional-probability/} \\
DP 1 & \nolinkurl{ib-dp-1-analysis-and-approaches-sl-sets-and-venn-diagrams/} \\
DP 1 & \nolinkurl{ib-dp-1-analysis-and-approaches-sl-probability/} \\
DP 2 & \nolinkurl{ib-dp-2-analysis-and-approaches-sl-discrete-random-variables/} \\
\bottomrule
\end{tabularx}
\caption{Probability chapters of the IB courses (October 2026). The DP chapters exist in the four courses (analysis and approaches, applications and interpretation, SL and HL); the SL analysis-and-approaches addresses are given.}
\label{tab:addresses-ib}
\end{table}

\end{document}